\documentclass[12pt]{article}
\usepackage{mathrsfs}
\usepackage{amsfonts}
\usepackage{amsmath, amsthm, amsfonts, amssymb}
\theoremstyle{definition}

\numberwithin{equation}{section} \theoremstyle{remark}

\def\<{\langle}
\def\>{\rangle}

\def\ra{\rightarrow}
\def\p{\partial}
\def\a{\alpha}

\def\D{{\cal D}}

\def\l{{\lambda}}

\def\CC{{\bf C}}

\def\-{\overline}
\def\ale{\mathrel{\mathop{<}\limits_{\sim}}}

\def\h{\hbox}
\def\d{\delta}
\def\endpf{\hbox{\vrule height1.5ex width.5em}}

\def\d{\delta}
\def\d{\delta}
\def\b{\beta}
\def\a{\alpha}
\def\endpf{\hbox{\vrule height1.5ex width.5em}}
\def\CC{\bf C}

\def\M*{\wt{M^*}}

\def\-{\overline}
\def\ale{\mathrel{\mathop{<}\limits_{\sim}}}

\def\D{\Delta}

\def\h{\hbox}

\def\wt{\widetilde}
\def\ra{\rightarrow}

\def\d{\delta}

\def\endpf{\hbox{\vrule height1.5ex width.5em}}

\def\a{\alpha}
\def\d{\delta}

\def\D{\Delta}
\def\d{\delta}
\def\b{\beta}
\def\ord{\hbox{Ord}}
\def\a{\alpha}

\def\endpf{\hbox{\vrule height1.5ex width.5em}}

\def\beq{\begin{equation}}
\def\nneq{\end{equation}}
\def\beqn{\begin{eqnarray}}
\def\neqn{\end{eqnarray}}
\def\beqna{\begin{eqnarray*}}
\def\neqna{\end{eqnarray*}}
\def\bedis{\begin{displaymath}}
\def\nedis{\end{displaymath}}
\def\-{\overline}
\begin{document}

\bigskip

\title{\bf  A Bishop surface  with a  vanishing Bishop invariant}

\author{Xiaojun Huang\footnote{
Supported in part by NSF-0500626} \ \  and\  Wanke Yin}


\date{March, 2007 }

\vspace{3cm}
\maketitle


\centerline{\bf  Abstract}
 We derive a complete set of invariants for a formal Bishop surface near a point of complex tangent with a vanishing Bishop invariant
 under the action of formal  transformations. We prove that the
 modular space of Bishop surfaces with a vanishing Bishop invariant and with a fixed  Moser invariant $s<\infty$ is of infinite dimension. We also prove that the equivalence class of the germ of a generic real analytic Bishop surface near a complex tangent
with a vanishing Bishop invariant can not be determined by a finite part of the Taylor expansion of its defining equation.
This answers, in the negative, a  problem raised   by  J. Moser in 1985 after his joint work with Webster in 1983 and his own work in 1985.
Such a phenomenon is  strikingly different from the celebrated theory of Moser-Webster for elliptic Bishop surfaces with non-vanishing
Bishop invariants. We also show that a formal map between two real
analytic Bishop surfaces with the Bishop invariant $\lambda=0$ and with the Moser invariant $s\not = \infty$ is convergent.
Hence, two real analytic Bishop surfaces with $\lambda=0$ and $s<\infty$ are holomorphically equivalent if and only if they have the same formal normal form (up to a trivial rotation).
Notice that there are many  non-convergent formal
transformations between  Bishop surfaces with $\lambda=0$ and $s=\infty$.
Notice also that a generic formal map between two real analytic  hyperbolic Bishop surfaces is divergent as shown by Moser-Webster and Gong.
Hence, Bishop surfaces with  a vanishing Bishop invariant and $s\not = \infty$  behave  very differently, in this respect,  from  hyperbolic Bishop surfaces or elliptic Bishop surfaces with $\lambda=0 $ and $ s=\infty$.
We also show that a Bishop surface with  $\lambda=0$ and $s<\infty$ generically has a trivial
automorphism group and
 has  the largest possible automorphism group if and only if it is biholomorphic to the model surface $M_s=\{(z,w)\in {\mathbb C}^2:\ w=|z|^2+z^s+\-{z}^s\}$.  Notice that, by the Moser-Webster theorem, an elliptic Bishop surface with $\lambda\not = 0$, always has automorphic group ${\mathcal Z}_2$.
Hence, Bishop surfaces with  $\lambda=0$ and  $s\not = \infty$
 have the similar character as  that of strongly pseudoconvex  real hypersurfaces in the complex spaces of higher dimensions.
\bigskip\bigskip

\section{Introduction and statements of main results}
In this paper, we study  the precise holomorphic structure of a real analytic Bishop surface near  a
complex tangent point with the  Bishop invariant vanishing. A Bishop surface is a
  generically embedded real surface in the complex space of
 dimension two. Points on  a  Bishop surface are either
 totally real  or have   non-degenerate complex tangents.
The holomorphic structure near a totally real point
is trivial. Near
a point with a complex tangent, namely, a point
with a non-trivial complex tangent space of type $(1,0)$,  the consideration
could be much more subtle. The study of this problem was initiated by
 the celebrated paper of Bishop  in 1965 [Bis], where for a point $p$ on a Bishop surface  $M$  with a complex tangent, he
defined an invariant $\lambda$ now called the Bishop invariant. Bishop showed that
there is  a  holomorphic change of variables, that maps $p$ to $0$, such that $M$, near
$p=0$, is defined in the complex coordinates $(z,w)\in {\mathbb C}^2$
by
\begin{equation}
w=z\-{z}+\lambda (z^2+\-{z}^2)+o(|z|^2),
\label{eqn:000}
\end{equation}
where $\lambda\in [0,\infty]$. When $\lambda=\infty$, (\ref{eqn:000}) is understood as $w=z^2+\-{z}^2+o(|z|^2).$
It is now a standard terminology to call $p$ an elliptic, hyperbolic or parabolic point of $M$, according to whether $
\lambda\in [0,1/2)$, $\lambda \in (1/2,\infty)$ or $\lambda=1/2,\infty$, respectively.

 Bishop discovered  an important geometry associated with $M$ near an elliptic complex tangent $p$ by
proving the existence of a family
of holomorphic disks attached to $M$ shrinking down to $p$.
He also proposed several problems
concerning the uniqueness and regularity
 of the geometric object  obtained by taking the union of all
locally attached holomorphic disks. These problems,
including their higher dimensional cases, were completely answered  through
the combining efforts of many people. (See [KW1], [BG], [KW2], [MW], 
[Mos], [HK], [Hu3]; in particular, see [KW1], [MW],  [Hu3]).

Bishop invariant is a quadratic invariant, capturing the basic geometric character of the surface.
The celebrated work of Moser-Webster [MW] first investigated  the
 more subtle higher order invariants.
Different from Bishop's approach of using the attached holomorphic
disks, Moser-Webster's starting point is the existence of a more
dynamically oriented object: an intrinsic pair of involutions on the
complexification of the surface near a non-exceptional complex
tangent. Here, recall that the Bishop invariant is said to be
non-exceptional if $\lambda\not = 0, 1/2, \infty$ or if
$\lambda\nu^2-\nu +\lambda=0$ has no roots of unity in the variable
$\nu$. Moser-Webster  proved that, near a non-exceptional complex
tangent, $M$ can always be mapped, at least, by a formal
transformation to the  normal form defined in the complex
coordinates $(z,w=u+iv)\in {\mathbb C}^2$ by:
\begin{equation}
u=z\-{z}+(\lambda +\epsilon u^s)(z^2+\-{z}^2)\ , \ v=0\ ,\
\epsilon\in\{0,1,-1\}\ ,\ s\in {\mathbb Z}^+.
\end{equation}

Moser-Webster also provided a convergence proof of the above mentioned formal transformation in
the non-exceptional elliptic case: $0<\l<1/2$. However, the intriguing elliptic case with $\lambda=0$ has to be excluded from their theory.
Instead, Moser in [Mos] carried out a study for $\lambda=0$ from a more formal power series point of view.
Moser derived the following  formal pseudo-normal form for $M$ with
$\lambda=0$:
\begin{equation}
w=z\-{z}+z^s+\-{z}^s+2Re\{\sum_{j\ge s+1}a_{j}z^j\}.
\end{equation}
Here $s$ is the simplest higher order invariant of $M$ at a complex
tangent with a vanishing Bishop invariant, which we call the Moser
invariant. Moser showed that when $s=\infty$, $M$ is then
holomorphically equivalent to the quadric $M_{\infty}=\{(z,w)\in
{\mathbb C}^2: w=|z|^2\}$.

Moser's formal pseudo-normal form is still subject  to the simplification of
a very complicated infinitely dimensional group $aut_{0}(M_{\infty})$, the formal
self-transformation group of $M_{\infty}$. And it was left open  from the work of Moser [Mos] to derive any higher order invariant
other than $s$ from the Moser pseudo-normal form. At this point, we mention that $aut_{0}(M_{\infty})$
 contains many non-convergent elements. Based on this, Moser asked two basic problems
concerning a Bishop surface near a vanishing Bishop invariant in his
paper [Mos]. The first one is concerning the analyticity of the
geometric object formed by the attached disks up to the complex
tangent point. This was answered in the affirmative in [HK].
Hence, the work of [HK], together with that of
Moser-Webster [MW], shows that, as far as the analyticity of the
local hull of holomorphy is concerned,  all elliptic Bishop surfaces
are of the same character. The second problem that Moser asked is
concerning the higher order invariants. Notice that by the
Moser-Webster normal form, an analytic elliptic Bishop surface with
$\lambda\not =0$ is holomorphically equivalent to an algebraic one
and possesses at most two more higher order invariants. Moser asked
if $M$ with $\lambda=0$ is of the same  character as that for
elliptic surfaces with $\lambda\not =0$.
 Is the equivalence class of a Bishop surface with $\lambda=0$  determined  by an algebraic surface obtained by truncating  the Taylor expansion of its defining equation at a sufficiently higher order level?
Gong showed in [Gon2] that under the equivalence relation
of a   smaller class of transformation group, called  the group of holomorphic symplectic transformations,
$M$ with $\lambda=0$ does have an infinite set of invariants. However, under this equivalence relation, elliptic surfaces with non-vanishing invariants
also have infinitely many invariants.
Gong's  work later on (see, for example, [Gon2-3] [AG])  demonstrates that as far as many dynamical properties are concerned, exceptional or non-exceptional hyperbolic,  or even parabolic complex tangents are not much different from each other.

In this paper, we derive a  formal normal form for a Bishop surface near a vanishing Bishop invariant, by introducing a quite different weighting system. This
new weighting system  fits extremely  well in our setting
and may have applications in many other problems.
We will obtain a complete set of invariants under the action of the formal transformation group. We show, in particular, that the modular space for  Bishop surfaces with a vanishing Bishop invariant and with a fixed (finite) Moser invariant $s$ is an infinitely dimensional manifold in a
Fr\`echet space. This then immediately provides an answer, in the negative, to Moser's problem concerning the determination of a Bishop surface with
a vanishing Bishop invariant from a finite truncation of its Taylor expansion.
Furthermore,
it can also be combined with some already known arguments  to show that most  Bishop surfaces with $\lambda=0,\ s\not =\infty$ are not holomorphically equivalent to  algebraic surfaces.
Hence, one sees a striking difference of an elliptic Bishop surface with a vanishing Bishop invariant from elliptic Bishop surfaces with non-vanishing Bishop
invariants.
The general phenomenon that the infinite dimensionality of the modular space has the consequence that   any subclass formed by a countable union of finite dimensional spaces is  of the first category in the modular space  seems already clear  even to Poincar\'e [Po]. In the CR geometry category, we refer the reader to a  paper of Forstneric [For]
in which the infinite dimensionality of the modular space of generic CR manifolds is used to show
that CR manifolds holomorphically equivalent to algebraic ones form a very thin set among all real analytic CR manifolds. Similar to what Forstneric did in  [For], our argument to show the generic non-algebraicity from the infinite dimensionality of the modular space also uses the Baire category theorem. 

It is not clear to us if   the new normal form obtained in this paper for  a real analytic
Bishop surface with $\lambda=0,\ s<\infty$ is always  convergent. However,
we will show that if the formal normal form is convergent, then the
map transforming the surface to its normal form must be convergent
in case the Moser invariant $s\not =\infty$. Remark that there are
many non-convergent formal maps transforming  real analytic Bishop
surfaces with a vanishing Bishop invariant and with $s=\infty$ to
the model surface  $M_{\infty}$ defined before. (See [MW] [Mos]
[Hu2]). Hence, our convergence theorem   reveals a non-trivial role
that the Moser invariant has played   in the study of the precise
holomorphic structure of a Bishop surface with $\lambda=0$.
At this point, we would like to mention that  there are many other
different problems where one  also considers
the convergence of formal power series, though very  different methods and approaches  need to be employed in different settings.
To name a few,  we here  mention the papers of
Baouendi-Ebenfelt-Rothschild [BER][BMR][MMZ], Webster [We], Stolovitch [St]
and the references therein. In the research described in [BER][BMR][MMZ],  one
tries to understand the convergence  of formal CR maps between not
too degenerate real analytic CR manifolds. In [We] [Sto], one
encounters other type of convergence problems in the normalization
of real submanifolds in ${\Bbb C}^n$.

Our convergence argument  uses the Moser-Webster [MW] polarization, as in the non-vanishing Bishop invariant case treated by Moser-Webster.
However, different from the Moser-Webster situation, we do not have a
pair of involutions, which were the starting point of the
Moser-Wbetser theory.
The  main idea in the present paper for dealing with our  convergence problem
is to find a new surface hyperbolic geometry, by making use of the flattening theorem of Huang-Krantz [HK].


We next state our main results, in which we will use some terminology to be defined in the next section:

\bigskip
\textbf{Theorem 1.1:}\hspace{0.2cm} {\it Let $M$ be a formal  Bishop surface with an elliptic complex tangent at $0$,
whose  Bishop invariant $\l$ is $0$  and  whose
 Moser invariant $s$ is a finite integer  greater than two . Then
there exists a formal transformation,
$$(z',w')=F(z,w)=(\wt{f}(z,w),\wt{g}(z,w)),\ \ F(0,0)=(0,0)$$
 such that
in the $(z',w')$
coordinates, $M'=F(M)$ is represented near the origin  by a formal equation of the following normal form:
\begin{equation*}
w'=z'\bar{z'}+z'^s+\bar{z'}^s+\varphi(z')+\overline{\varphi(z')}
\end{equation*}
where
$$
\varphi(z')=\sum\limits_{k=1}^{\infty}\sum\limits_{j=2}^{s-1}a_{ks+j}z'^{ks+j}.
$$
Such a formal transform is unique up to a composition from the left with a
rotation  of the form:
$$z''=e^{i\theta}z',\ w''=w',\ \ \hbox {where } \theta \ \hbox{is a constant with }\ e^{is\theta}=1.$$
}
\bigskip

\textbf{Theorem 1.2}:\hspace{0.2cm} {\it Let $M$ and $M'$ be  real analytic Bishop surfaces near $0$
 with  the Bishop invariant vanishing and the Moser invariant finite.
Suppose that
$F:(M,0) \longrightarrow (M',0)$ is a formal equivalence map. Then F is biholomorphic
near $0$. }
\bigskip

Define ${\mathcal Z}_s$ for the group of transformations consisting of maps of
the form $\{\psi_\theta: (z,w) \mapsto (e^{i\theta}z,w), \hskip 5pt
e^{is\theta}=1\}$.

We next give  several immediate  consequences of Theorems 1.1 and 1.2:
\bigskip

\textbf{Corollary 1.3:}\hspace{0.2cm} {\it (a): Suppose $M_{nor}$ is a formal Bishop surface near the origin defined by
$${w=z\bar{z}+z^s+\bar{z}^s+2Re\{\sum\limits_{k=1}^{\infty}\sum\limits_{j=2}^{s-1}
a_{ks+j}z^{ks+j}}\}.
$$
Then the group of  the origin preserving formal self-transformations of $M_{nor}$, denoted by $\hbox{aut}_0(M_{nor})$,  is a subgroup of
 ${\mathcal  Z}_s$. Moreover,  $\psi_\theta \in \hbox{aut}_0(M_{nor})$ if and
only if
$$a_{ks+j}=0\ \ \hbox{for any k and j with } k\ge 1,\ 2\le j\le s-1,\ \ e^{\sqrt{-1}j\theta}\not =1.$$
\medskip
\noindent
(b):\  $\hbox{aut}_0(M_s)={\mathcal Z}_s$, where $M_s$ is defined by $w=z\-{z}+z^s+\-{z}^s$.\\
\medskip
\noindent
(c):\ Any subgroup of ${\mathcal Z}_s$ can be realized  as the formal automorphism group of
a certain $M_{nor}$. \\
\medskip
\noindent
(d):\ Let $M$ be a formal Bishop surface   with a vanishing Bishop invariant and $s<\infty$ at $0$. Then
$\hbox{aut}_0(M)$ is isomorphic  to a subgroup of ${\mathcal Z}_s$.
\medskip

\noindent
(e):\  Let $M$ be a real analytic Bishop surface with a vanishing Bishop invariant and the Moser invariant $s<\infty$ at $0$. Suppose that $aut_{0}(M)={\mathcal Z}_s$. Then $(M,0)$ is biholomorphic to $(M_s,0)$, where $M_s$, as before, is defined by
$w=z\-{z}+z^s+\-{z}^s$.

\noindent
(f):\  Let $M$ be a real analytic elliptic Bishop surface with  $\l=0 $ and $s$ a prime number at $0$. Then $aut_{0}(M)$ is a trivial group  unless $(M,0)$ is biholomorphic to $(M_s,0)$.}
\bigskip
\bigskip

{\bf Corollary 1.4}: {\it Let $M_1$ and $M_2$ be  real analytic Bishop surfaces with $\lambda=0$ and $s\not =\infty$ at $0$. Suppose that $M_1$ has
a formal normal form:
$$w'=z'\bar{z'}+{z'}^s+\bar{z'}^s+2Re\{\sum\limits_{k=1}^{\infty}\sum\limits_{j=2}^{s-1}
a_{ks+j}{z'}^{ks+j}\};
$$
and suppose that $M_2$ has a formal normal form:
$$w'=z'\bar{z'}+{z'}^s+\bar{z'}^s+2Re\{\sum\limits_{k=1}^{\infty}\sum\limits_{j=2}^{s-1}
b_{ks+j}{z'}^{ks+j}\}.$$
 Then $(M_1,0)$ is biholomorphic to $(M_2,0)$ if and only if there is a constant $\theta$, with $e^{s\theta\sqrt{-1}}=1$,
such that $a_{ks+j}=e^{\theta j\sqrt{-1}}b_{ks+j}$ for any
$k\ge 1$ and $j=2,\cdots,s-1$.}
\bigskip

{\bf Theorem 1.5}: {\it A generic real analytic Bishop surface with
a vanishing Bishop invariant and $s\not =\infty$ is not
holomorphically equivalent to an algebraic surface in
${\mathbb{C}}^2$.}

\bigskip

{\bf Acknowledgment}: The key part of this work was completed when the first author was visiting, in
January of 2006, the School of Mathematics, Wuhan University, China and  when both authors were enjoying the month long visit at the Institute of Mathematical Sciences, The Chinese University of Hong Kong in the Spring of 2006.
The first author would like  very much to thank his friends Professors Hua Chen and Gengsheng Wang at Wuhan University
for their hospitality  during the visit. Both authors would also like to  express their appreciation to IMS at the Chinese University of Hong Kong
for its generous  supports and helps provided during the authors' visit.

\bigskip

\section{Uniqueness of  formal maps between approximately   normalized surfaces }

In what follows, we use $(z,w)$ or $(z',w')$ for the coordinates for ${\mathbb C}^2$.
Let $A(z,\-{z})$ be a formal power series in $(z,\-{z})$ without constant term. We say that the order of $A(z,\bar{z})$ is $k$ if $A(z,\-{z})=\sum_{j+l=k}A_{j\-{l}}z^j\-{z}^l+o(|z|^k)$ with at least one of the $A_{j\-{l}}\in {\mathbb C}$ ($j+l=k)$ not equal to $0$. In this case,
we write  $\ord(A(z,\-{z}))=k$. We say $\ord(A(z,\-{z}))\ge k$ if $A(z,\-{z})=O(|z|^k)$.

Consider a formal real surface M in $\mathbb{C}^2$ near the origin. Suppose that $0$ is a point of complex tangent for $M$. Then, after a linear change of
variables, we can assume that $T^{(1,0)}_0M=\{w=0\}$. If there is no change of coordinates such that $M$ is defined by an equation of the form $w=O(|z|^3)$, we then say $0$ is a point of $M$ with a non-degenerate complex tangent. In this case, Bishop showed that
there is a change of coordinates in which $M$ is defined by ([Bis] [Hu1])
\begin{equation}
w=z\-{z}+\lambda (z^2+\-{z}^2)+O(|z|^3).
\end{equation}
Here $\lambda\in [0,\infty]$ and when $\lambda=\infty$, the equation takes the form: $w=z^2+\-{z}^2+O(|z|^3).$
$\lambda$ is the first absolute invariant of $M$ at $0$, called the Bishop invariant. Bishop invariant is a quadratic invariant, resembling to
the Levi eigenvalue in the hypersurface case. When $\lambda\in [0,1/2)$, we say that $M$ has an elliptic complex tangent at $0$. In this paper, we are only interested in the case  of an elliptic complex tangent.
We need only to
study  the case of $\lambda=0$; for, in the case with $\lambda\in (0,1/2)$, the surface has been well understood by the work of Moser-Webster [MW].
When $\lambda=0$, Moser-Webster and Moser   showed in [MW] [Mos] that there is an integer $s\ge 3$ or $s=\infty$ such that $M$ is defined by
\begin{equation}
w=z\bar{z}+z^s+\bar{z}^s+E(z,\bar{z}),
\label{eqn:Jam000}
\end{equation}
where    $E$ is a formal power series in $(z,\-{z})$  with $\ord(E) \geq s+1.$
When $s=\infty$, we understand the defining equation as $w=z\-{z}$, namely, $M$ is formally equivalent to the quadric $M_{\infty}=\{w=z\-{z}\}$.
$s$ is the next absolute invariant for $M$, called the Moser invariant.
The case for $s=\infty$ is also well-understood through the work of Moser [Mos]. Hence, in all that follows, our $M$ will have $\lambda=0$ and
a fixed $s<\infty$.


A  formal map
$z'=F(z,w),\ w'=G(z,w)$ without constant terms is called
an invertible  formal transformation (or simply, a formal transformation) if $\frac{\partial (F,G)}{\partial
(z,w)}(0,0)$ is invertible. When a formal map has no constant term, we also  say that it preserves the origin.


\bigskip

\textbf{Lemma 2.1}:\hspace{0.2cm} {\it Let $M$ be defined as in (\ref{eqn:Jam000}). Suppose that $z'=F(z,w),w'=G(z,w)$ is a
formal transformation preserving the origin and  sending $M$ into $M'$, which is  defined by
$w'=z'\-{z}'+E^*(z',\-{z'})$. Then

\noindent (i): $F=az+bw+O(|(z,w)|^2),\ G=cw+O(|w|^2+|zw|+|z|^3)$
where $c=|a|^2, a \neq 0$.

\noindent
(ii): Suppose  that $M$ and $M'$ are further defined by  $w=E(z,\-{z})=z\-{z}+z^s+\-{z}^{s}+o(|z|^s)$ and
$w'=E^{*}(z',\bar{z'})=z'\-{z'}+{z'}^s+\-{z'}^{s}+o(|z'|^s)$, respectively. Here $s\ge 3$. Then
$$F=(e^{i\theta}z+O(|z|^2+|w|),w+O(|w|^2+|zw|+|z|^3)),\ \hbox{ where}\ \theta \ \hbox{is a constant with }\ e^{is\theta}=1 .$$

\noindent
(iii): In (i), when  $\-{E(z,\-{z})}=E(z,\-{z})$ and  $\overline{E^*(z',\-{z}')}=E^*(z',\-{z}')$, we have
$G(z,w)=G(0,w)$.}

\bigskip

{\it Proof of Lemma 2.1}: (i) is the  content  of Lemma 3.2 of
[Hu1]. To prove (ii), we write
$F=(az+f,cw+g)$, where
by (i), we can assume that
$$
f(z,w)=O(|z|^2+|w|)\ \ ,\ \
g(z,w)=O(|w|^2+|zw|+|z|^3).
$$
Notice that
$$
f(0,E(0,\bar{z}))=O(\bar{z}^s)\ \ ,\ \
\overline{{f}}(\bar{z},\bar{E}(\bar{z},0))=O(\bar{z}^2)\ \
,\ \ g(0,E(0,\bar{z}))=o(\bar{z}^s).
$$
Applying  the defining equation of $M'$ , we have, on $M$, the following:
\begin{equation*}
\begin{array}{lll}
cw+g(z,w)&=&|a|^2|z|^2+\bar{a}\bar{z}f(z,w)
   +az\overline{f}(\bar{z},\bar{w})
   +f(z,w)\overline{f}(\bar{z},\bar{w})\\
&&+\left(az+f(z,w)\right)^s+\left(\bar{a}\bar{z}
   +\overline{f}(\bar{z},\bar{w})\right)^s+o(|z|^s).
\end{array}
\end{equation*}
Regarding $z$ and $\bar{z}$ as independent variables in the above equation and then letting $z=0,w=E(0,\bar{z}),\bar{w}=\bar{E}(\bar{z},0)$,
 we obtain
$$
c\bar{z}^s+o(\bar{z}^s)=(\bar{a}\bar{z})^s+o(\bar{z}^s).
$$
Hence, it follows that $c=\bar{a}^s$. Together with $c=|a|^2$ and $s
\geq 3$, we get
$$
c=1\ ,\ a=e^{i\theta}\ ,\ \hbox{where $\theta$ is a constant}.
$$
Now we turn to the proof of (iii). Notice that
$$
G(z,w)=|F(z,w)|^2+E^*(F(z,w),\overline{F(z,w)})\hskip 5pt for \hskip
5pt (z,w) \in M.
$$
Since $E^*$ is now assumed to be  formally real valued, we have
$$
G(z,w)=\overline{G(z,w)} \ \ \   \hbox{on}\   M.
$$
Write
$$
G(z,w)=\sum\limits_{\a,\b}^{\infty}a_{\a\b}z^\a w^\b.
$$
We will prove  inductively that $a_{\a\b}=\-{a_{\a\b}}$ for $\a=0$ and $a_{\a\b}=0$ otherwise.
First, for each $m>>1$, write $E=E_{(m)}(z,\-{z})+E_m$ with $E_{(m)}(z,\-{z})$ a polynomial of degree at most $m-1$ and $E_m=O(|z|^m)$.
Then for any $m>>1$, there are integers $N_1(m)>>m$ and $N_2(m)>>m$ such that
\begin{equation}
\sum\limits_{\a,\b=0}^{N_2(m)}a_{\a\b}z^\a
w^\b=\sum\limits_{\a,\b=0}^{N_2(m)}\-{a_{\a\b}z^\a}w^\b +o(|z|^m)\ \
,\ \ w=z\-{z}+E_{(N_1(m))}(z,\-{z}).
\end{equation}
Next, suppose that $N_0=\a_0+2\beta_0$ is the smallest number
 such that $a_{\a\b}$ is real-valued for $\a=0$, and  zero otherwise whenever $\a+2\b<N_0$. (If such an $N_0$ does not exist, then Lemma 2.1 (iii) holds automatically).
Choose $m>>N_0$. For $0<r<<1$, define $\sigma_{N_1}(\xi,r)$ to be the biholomorphic map from the unit disk in $\mathbb C$ to the smoothly bounded simply connected domain:
$\{\xi\in {\mathbb C}:\ |\xi|^2+r^{-2}E_{(N_1)}(r\xi,r\-{\xi})<1\}$ with $\sigma_{N_1}(\xi,r)=\xi(1+O(r)).$
Since the disk $(r\sigma_{N_1}(\xi,r),r^2)$ is attached to $M_{N_1}$ defined by $ w=z\-{z}+E_{(N_1)}(z,\-{z})$,
it follows that
\begin{equation}
\sum\limits_{\a+2\b=N_0}a_{\a\b}r^{N_0}\xi^\a =\sum\limits_{\a+2\b=N_0}\-{a_{\a\b}\xi^\a}r^{N_0}+o(r^{N_0}), \ \ |\xi|=1.
\end{equation}
Letting $r\ra 0$, we get
\begin{equation}
\sum\limits_{\a+2\b=N_0}a_{\a\b}\xi^\a =\sum\limits_{\a+2\b=N_0}\-{a_{\a\b}\xi^\a}, \ \ |\xi|=1,
\end{equation}
from which we see that when $\a+2\b=N_0$, $a_{\a\b}$ is real for $\a=0$, and zero otherwise. This contradicts the choice of $N_0$ and thus completes the proof of Lemma 2.1 (iii).
$\endpf$

\bigskip
The main purpose of this section is to prove the following uniqueness result for mappings between approximately normalized surfaces:
\bigskip

\textbf{Theorem 2.2}:\hspace{0.2cm} {\it Suppose that the formal power series
\begin{equation}
\left\{
\begin{array}{ll}
z'=z+f(z,w),\hskip 10pt&f(z,w)=O(|w|+|z|^2) \\
w'=w+g(w), &g(w)=O(|w|^2)
\end{array}
\right.
\end{equation}
transforms the formal Bishop surface $M$ defined by
$$
w=z\bar{z}+2Re\left(z^s+\sum\limits_{k=1}^{n}\sum\limits_{j=2}^{s-1}a_{ks+j}z^{ks+j}\right)+E_1(z,\-{z})
$$
to the formal Bishop surface defined by
$$
w'=z'\bar{z'}+2Re\left(z'^s+\sum\limits_{k=1}^{n}\sum\limits_{j=2}^{s-1}b_{ks+j}z'^{ks+j}\right)+E_2(z',\-{z'})
$$
where $n\ge 1$, $a_{ks+j},b_{ks+j}$ are complex numbers, and $E_1(z,\-{z}),\ E_2(z,\-{z})=o(|z|^{ns+s-1}).$
Then $f(tz,t^2w)=O(t^{2n+1})$, $g(t^2w)=O(t^{2n+2})$, as $t\in {\mathbf R}\ra 0$, and
$a_{ks+j}=b_{ks+j}$ for all $k\le n$ and $j=2,\cdots, s-1$.}
\bigskip

One of the crucial ideas for the proof of Theorem 2.2 is to set the weight of $\-{z}$ differently from that of $z$. More precisely,
we set the weight of $z$ to be $1$ and that of $\bar{z}$ to be $s-1$. For a formal power series  $A(z,\-{z})$ with no constant term,
we say that $wt(A(z,\-{z}))=k$, or $wt(A(z,\-{z}))\ge k$, if $A(tz,t^{s-1}\-{z})=t^kA(z,\-{z})$, or , $A(tz,t^{s-1}\-{z})=O(t^k)$,
respectively, as $t\in {\mathbb R}\ra 0$.
 In all that follows, we  use
$\Theta_l^j$ to denote  a formal power series in $z$ and $\bar{z}$ of
order at least  $j$ and weight at least $l$. (Namely, $\Theta^j_l(tz,t\-{z})=O(t^j)$ and
 $\Theta_l^j(tz,t^{s-1}\-{z})=O(t^l)$ as $t\ra 0$).
We use $\mathbb{P}_l^j$ to denote  a homogeneous polynomial in $z$ and $\bar{z}$ with the exact order $j$ and
weight at least $l$. {\it We emphasize that $\Theta_l^j$ and $\mathbb{P}_l^j$ may be different in  different contexts}.

In what follows, we also define the normal weight of $z, w$ to be $1,2$, respectively.
For a formal power series $h(z,w,\-{z},\-{w})$, we use $wt_{nor}(h)\ge k$ to denote the vanishing property:
$h(tz,t^2w,t\-{z},t^2\-{w})=O(t^k)$ as $t\ra 0$.
Let $h(z,w)$ be a formal power series in $(z,w)$ without constant term.
Then we have the formal expansion:
$$
h(z,w)=\sum\limits_{l=1}^{\infty}h_{nor}^{(l)}(z,w)
$$
where
$$
h_{nor}^{(l)}(tz,t^2w)=t^lh_{nor}^{(l)}(z,w)
$$
is a polynomial in $(z,w)$.
Notice that $h_{nor}^{(l)}(z,w)$ is  homogeneous of degree $l$ in the standard weighting system which assigns the weight of
$z$ and $w$ to be $1$ and $2$, respectively.
In what follows, we write
\begin{equation}
h_l(z,w)=\sum\limits_{j=l}^{\infty}h_{nor}^{(j)}(z,w)\  \hbox{ and
}\ h_{(l)}=\sum_{j=1}^{l-1}h_{nor}^{(j)}(z,w).
\label{eqn:James-00010}
\end{equation}

\bigskip
{\it Proof of Theorem 2.2}:
We need to prove that any solution $(f,g)$ of the following equation
has
 the property that
 $wt_{nor}(f(z,w))\ge 2n+1$, $wt_{nor}(g(w))\ge 2n+2$
under the normalization conditions as in the theorem:
\begin{equation}
\begin{array}{lll}
w+g(w)&=&(z+f(z,w))(\bar{z}+\overline{f(z,w)}) +2Re \large\{
(z+f(z,w))^s\\
&&+\sum\limits_{k=1}^{n}\sum\limits_{j=2}^{s-1}b_{ks+j}(z+f(z,w))^{ks+j}
\large\}+E_2(f(z,w),\-{f(z,w)})
\end{array}
\label{eqn:Jam01}
\end{equation}
where  $w=z\bar{z}+z^s+\bar{z}^s+E(z,\bar{z})$ with
$$E=2Re\left(\sum\limits_{k=1}^{n}\sum\limits_{j=2}^{s-1}a_{ks+j}z^{ks+j}\right)+E_1(z,\-{z}).$$
With an immediate simplification, (\ref{eqn:Jam01}) takes the form:
\begin{equation}
\begin{array}{lll}
g(w)&=&\bar{z}f(z,w)+z\overline{f(z,w)}+|f(z,w)|^2
+2Re\big\{(z+f(z,w))^s-z^s\\
&&+\sum\limits_{k=1}^{n}\sum\limits_{j=2}^{s-1}\left(b_{ks+j}(z+f(z,w))^{ks+j}-a_{ks+j}z^{ks+j}\right)\big\}
+o(|z|^{ns+s-1})
\end{array}
\label{eqn:Jam02}
\end{equation}
\medskip

In the proof of Theorem 2.2, we set the following convention. For any positive integer $N$, we define $a_N$ and $b_{N}$ to be as in Theorem 2.2 if
$N=ks+j$ with $k\le n,\ 2\le j\le s-1$, and to be $0$ otherwise.
For the rest of this section, we  will define a positive integer $N_0$ as follows:

Suppose that there is a pair of integers $(j_0,k_0)$
 such that $s<k_0s+j_0(\le ns+s-1)$ is the smallest
 number satisfying $a_{k_0s+j_0} \neq b_{k_0s+j_0}.$ We then define $N_0=k_0s+j_0$.
 Otherwise, we  define $N_0=sn+s$.

The proof of Theorem 2.2 is carried out in two steps, according
to  the vanishing order of $f$  being even or odd.
\medskip

\textbf{Step I of the proof of Theorem 2.2:} \ \
In this step, we assume that either
$$\ord\left(f(z,w(z,\bar{z})\right)=2t$$ is an even number or $f\equiv 0$,  where
$w(z,\bar{z})=z\bar{z}+z^s+\bar{z}^s+E(z,\bar{z})$.
Write $g(w)=c_lw^l+o(w^l)$.

Denote by $\widehat{N_0}=\hbox{min}\{N_0,\ \ord(f), sn+s-1\}$. (If $f\equiv
0$, we define $\ord(f)=\infty$.)
 Then (\ref{eqn:Jam02}) gives the following:
\begin{equation}
c_lz^l\bar{z}^l+O(|z|^{2l+1})=2Re[(b_{N_0}-a_{N_0})z^{N_0}]+O(|z|^{\widehat{N_0}+1}).
\end{equation}
From this, we can easily conclude the following:

\medskip
\noindent
(2.I).  Suppose that  $2t\ge N_0$ and $c_l\not =0$. Then $2l>\min\{N_0,sn+s-1\}$  and $b_{N_0}=a_{N_0}$. By our choice of $N_0$, $N_0$ must be $ns+s$. Hence, the theorem in this case
readily follows.

\medskip
\noindent
(2.II). When $2t<N_0$, then $2l\ge \min\{2t+2, sn+s\}$ under the assumption that $c_l\not = 0$. Thus $l>t\ge 1$ (if $c_l\not = 0$).

Suppose that  $N_0=2t+1$ in Case (2.II). Assuming that $N_0<ns+s$ and   collecting terms with degree $2t+1$ in $(\ref{eqn:Jam02})$, we obtain

\begin{equation}
\bar{z}f_{nor}^{(2t)}(z,z\bar{z})+z\overline{f_{nor}^{(2t)}(z,z\bar{z})}+2Re\left((b_{N_0}-a_{N_0})z^{N_0}\right)=0
\end{equation}
This clearly forces that $a_{N_0}=b_{N_0}$. Thus, we must have $N_0=ns+s$ and  Theorem 2.2 also follows easily in this setting.  {\it Hence, we will assume, in what follows:
\medskip

\noindent (2.III). $ns+s>N_0\ge 2t+2,\ l>t\ge 1$.}
\medskip



Collecting   terms
with (the ordinary) degree 2t+1 in (\ref{eqn:Jam02}), we get:
\begin{equation}
\bar{z}f_{nor}^{(2t)}(z,z\bar{z})+z\overline{f_{nor}^{(2t)}(z,z\bar{z})}=0
\label{eqn:Jam03}
\end{equation}
Writing $f_{nor}^{(2t)}(z,w)=\sum\limits_{k+2l=2t}a_{kl}z^kw^l$ and
substituting it back to (2.12), we then get:\\
$$
f_{nor}^{(2t)}(z,w)=aw^t-\bar{a}z^2w^{t-1}
$$
for  $a \neq 0$. Hence
\begin{equation}
f(z,w)=f_{nor}^{(2t)}(z,w)+f_{2t+1}(z,w)=aw^t-\bar{a}z^2w^{t-1}+f_{2t+1}(z,w)
\label{eqn:Jam03-00}
\end{equation}
Next, a simple computation shows that
$wt(w) \geq s,\ \ord(w(z,\-{z})) \geq 2, \ wt(f_{nor}^{(2t)}) \geq
st+2-s,\ wt(\overline{f_{nor}^{(2t)}}) \geq
st,\ g=g_{2t+2},\ f=f_{nor}^{(2t)}+f_{2t+1}(z,w).$  Also   if $ \ l_1+l_2 \geq s$ with $\ l_2>
1,$ or $ l_1+l_2 >s$ with $l_2 \ge 1$, then $wt(z^{l_1}f_{nor}^{{(2t)l_2}})=l_1+l_2(ts+2-s) \geq ts+2$.
Moreover, $wt(z^{l_1}f_{nor}^{{(2t)}l_2}f_{2t+1}^{l_3}) \geq s$
if $l_1+l_2+l_3 \geq s-1,\ l_2^2+l_3^2 \neq 0$.

We can verify the following
$$|f(z,w)|^2=2Re(\-{f_{nor}^{(2t)}}f_{2t+1})+\Theta^{2t+2}_{st+2}+\Theta^2_{st+2}f_{2t+1}.$$

Substituting (\ref{eqn:Jam03-00}) into (\ref{eqn:Jam02}), we get:
\begin{equation}
\begin{array}{rcl}
g_{2t+2}(w)&=&2Re\{(\-{z}+sz^{s-1})f\}+|f(z,w)|^2+
            2Re\{\sum_{l=2}^{s}
            \mathbb{P}_{s-l}^{s-l}f^{l}\}\\
&&+2Re(\sum\limits_{
{ \tau=ks+j<
            N_0}}\sum_{l=0}^{\tau-1}
            \mathbb{P}_{l}^{l}f^{\tau-l})
             +2Re\left((b_{N_0}-a_{N_0})z^{N_0}\right)+\Theta^{N_0+1}_{N_0+1}\\
           &=&2Re\{(\bar{z}+sz^{s-1})f_{nor}^{(2t)}+(\bar{z}+sz^{s-1}
+\-{f_{nor}^{(2t)}})f_{2t+1}(z,w)\} \\ &
&+2Re\left((b_{N_0}-a_{N_0})z^{N_0}\right)
+\Theta_{s}^{2}f_{2t+1}(z,w)+\Theta_{s}^{2}\-{f_{2t+1}(z,w)}
+\Theta_{N_s}^{2t+2}
\end{array}
\end{equation}
Here  $N_0$ is defined as before and $N_s:= \hbox{min}\{ts+2, N_0+1\}$.

Notice that
\begin{equation}
\begin{array}{l}
\bar{z}f_{nor}^{(2t)}+z\overline{f_{nor}^{(2t)}}+2Re\{sz^{s-1}f_{nor}^{(2t)}\}\\
=2Re\{\bar{z}(aw^t-\bar{a}z^2w^{t-1})+sz^{s-1}(aw^t-\bar{a}z^2w^{t-1})\}\\
=-\bar{a}z^2\bar{z}w^{t-1}+z\bar{a}w^t-sz^{s-1}\bar{a}z^2w^{t-1}+\Theta_{ts+2}^{2t+2}\\
=(1-s)\bar{a}z^{s+1}w^{t-1}+\Theta_{ts+2}^{2t+2}
\end{array}
\end{equation}

Hence, we obtain
\begin{equation}
\begin{array}{rll}
g_{2t+2}(w)&=&(1-s)\bar{a}z^{s+1}(z\bar{z}+z^s)^{t-1}+(\bar{z}+sz^{s-1}+\Theta_{s}^{2})f_{2t+1}(z,w)\\
&&+2Re\left((b_{N_0}-a_{N_0})z^{N_0}\right)
            +(z+s\bar{z}^{s-1}+\Theta_s^2)\overline{f_{2t+1}(z,w)}\\
&&+2Re\{\-{f_{nor}^{(2t)}}f_{2t+1}(z,w)\} + \Theta_{N_s}^{2t+2}
\label{eqn:Jam05}
\end{array}
\end{equation}
If $t=1$, collecting terms of degree  $s+1$ in (\ref{eqn:Jam05}) and noticing that $N_0>s+1$ by the given condition, we get
\begin{equation}
\begin{array}{lll}
\sum_{2j}\delta_{2j}^{s+1}g_{nor}^{(2j)}(z\-{z})&=&
(1-s)\bar{a}z^{s+1}+\bar{z}f_{nor}^{(s)}(z,z\-{z})+z\overline{f_{nor}^{(s)}(z,z\-{z})}\\
&&
+az^2\-{f_{nor}^{(s-1)}(z,z\-{z})}+\-{az^2}f_{nor}^{(s-1)}(z,z\-{z})+\mathbb{P}_{s+2}^{s+1}.
\end{array}
\end{equation}
Here $\d_{2j}^{s+1}$ takes value $1$, when $2j=s+1$, and $0$ otherwise.

Since $s+2 \geq s+1$, $\mathbb{P}_{s+2}^{s+1}=\bar{z}A$ with $A$ a polynomial. Thus it follows easily  that
$(1-s)\bar{a}z^{s+1}$ divides  $\bar{z}$. This is a contradiction and thus $t
>1$. In particular, (\ref{eqn:Jam05}) can be written as

\begin{equation}
\begin{array}{lll}
g_{2t+2}(w)&=&(1-s)\bar{a}z^{s+1}(z\bar{z}+z^s)^{t-1}+(\bar{z}+sz^{s-1}+\Theta_{s}^{2})f_{2t+1}(z,w)\\
&&+2Re\left((b_{N_0}-a_{N_0})z^{N_0}\right)
            +(z+s\bar{z}^{s-1}+\Theta_s^2)\overline{f_{2t+1}(z,w)}+
\Theta_{N_s}^{2t+2}
\label{eqn:Jam06}
\end{array}
\end{equation}

We next prove the following:
\bigskip

\textbf{Lemma 2.3}:\hspace{0.2cm} Assume that $2t+j(s-2)+2 \leq m \leq
2t+(j+1)(s-2)+1$ with $0 \leq j \leq t-1$ and  $m\le N_0$. Then
\begin{eqnarray}
\begin{array}{lll}
g_m(w)&=&\bar{a}(1-s)^{j+1}z^{(j+1)s+1}(z\bar{z}+z^s)^{t-j-1}+(\bar{z}+sz^{s-1}
          +\Theta_{s}^{2})f_{m-1}(z,w)\\
       &&+(z+s\bar{z}^{s-1}+\Theta_s^2)\overline{f_{m-1}(z,w)}+2Re\left((b_{N_0}-a_{N_0})z^{N_0}\right)+\Theta_{N_s}^{m}
\label{eqn:Jam060}
\end{array}
\end{eqnarray}
\bigskip

{\it Proof of Lemma 2.3}: The argument  presented above gives   the proof of the lemma with $m=2t+2$.
    We complete the proof of the lemma in three steps.
\medskip

{\bf Step I of the proof of Lemma 2.3}:\ This step is not needed when $s=3$.
 Denote  $ m_0=2t+j(s-2)+2$, where $j$ is an integer with
$0\le j\le t-1$. Suppose that $m_0\le N_0$.
We also assume that there is an integer $m$ such that $m\ge m_0$, $m+1\le 2t+(j+1)(s-2)+1$ (such an $m$ certainly does not exist
if $s=3$), $m+1\le N_0$ and moreover the formula  (\ref{eqn:Jam060}) holds for this $m$.
Collecting terms of degree $m$ in (\ref{eqn:Jam060}), we get
\begin{equation}
g^{(m)}(z\bar{z})=\bar{z}f_{nor}^{(m-1)}(z,z\bar{z})+z\overline{f_{nor}^{(m-1)}(z,z\bar{z})}+\hat{\mathbb{P}}_{N_s}^{m}
\label{eqn:Jam01-1}
\end{equation}
Notice that $\hat{\mathbb{P}}_{N_s}^{m}(=\mathbb{P}_{N_s}^{m})$ must be real valued, and notice that
$g^{(m)}(z\-{z})$ is also of weight at least $N_s$.
We can write
\begin{equation}
g^{(m)}(z\bar{z})-\mathbb{P}_{ts+2}^{m}=\sum\limits_{\stackrel{\alpha+\beta=m}{\alpha+\beta(s-1)\geq
N_s}} a_{\alpha\bar{\beta}}z^{\alpha}\bar{z}^{\beta}
\label{eqn:James01}
\end{equation}
Write\\
\begin{equation}
\begin{array}{ll}
f_{nor}^{(m-1)}(z,z\-{z})=\sum\limits_{\widetilde{\alpha}+2\widetilde{\beta}=m-1}
         b_{\widetilde{\alpha}\widetilde{\beta}}z^{\widetilde{\alpha}}(z\bar{z})^{\widetilde{\beta}}
         =\sum\limits_{\widetilde{\alpha}+2\widetilde{\beta}=m-1}
         b_{\widetilde{\alpha}\widetilde{\beta}}z^{\widetilde{\alpha}+\widetilde{\beta}}\bar{z}^{\widetilde{\beta}}.
\label{eqn:Jam07}
\end{array}
\end{equation}
Then
\begin{equation}
\sum\limits_{\widetilde{\alpha}+2\widetilde{\beta}=m-1}
b_{\widetilde{\alpha}\widetilde{\beta}}z^{\widetilde{\alpha}+\widetilde{\beta}}
\bar{z}^{\widetilde{\beta}+1}+\sum\limits_{\widetilde{\alpha}+2\widetilde{\beta}=m-1}
\overline{b_{\widetilde{\alpha}\widetilde{\beta}}}\bar{z}^{\widetilde{\alpha}+\widetilde{\beta}}
z^{\widetilde{\beta}+1} =
\sum\limits_{\stackrel{\alpha+\beta=m}{\alpha+\beta(s-1)\geq N_s}}
  a_{\alpha\bar{\beta}}z^{\alpha}\bar{z}^{\beta}
\end{equation}
We see that if $m$ is
even, then $2b_{\widetilde{\alpha}\widetilde{\beta}}=a_{\alpha\bar{\beta}}+ic$
when $\alpha=\beta={m \over
2},\ \widetilde{\alpha}=1,\ \widetilde{\beta}={m \over 2}-1,\ c \in
\mathbb{R}$. The other relations are as follows:
\begin{equation}
b_{\widetilde{\alpha}\widetilde{\beta}}=b_{\alpha\bar{\beta}},\ \ \hbox{if}\
\widetilde{\alpha}+\widetilde{\beta}=\alpha,\
\widetilde{\alpha}+2\widetilde{\beta}=m-1,\ \widetilde{\beta}+1=\beta,\widetilde{\alpha}
> 1,\ \alpha+(s-1)\beta \geq N_s.
\end{equation}
From this, one can easily see that
\begin{equation}
wt(f_{nor}^{(m-1)}(z,\-{z}))\ge \hbox{min}\{\wt{\a}+\wt{\b}+(s-1)\wt{\b}\}=\hbox{min}\{{\a}+(s-1){\b}-s+1\}\ge N_s-s+1.
\label{eqn:Jam08}
\end{equation}

Substituting (\ref{eqn:Jam07}) into (\ref{eqn:Jam060}), we get
\begin{equation}
\begin{array}{lll}
g_{m+1}(w)&=&(1-s)^{j+1}\bar{a}z^{(j+1)s+1}(z\bar{z}+z^s)^{t-j-1}+(\bar{z}+sz^{s-1}
          +\Theta_{s}^{2})f_m(z,w)\\
       &&+(z+s\bar{z}^{s-1}+\Theta_s^2)\overline{f_m(z,w)}+\Theta_{N_s}^{m+1}
       +(sz^{s-1}+\Theta_{s}^{2})f_{nor}^{(m-1)}\\
       &&+2Re\left((b_{N_0}-a_{N_0})z^{N_0}\right)+(s\bar{z}^{s-1}+\Theta_{s}^{2})\overline{f_{nor}^{(m-1)}}
\end{array}
\end{equation}
By (\ref{eqn:Jam08}), we get
$$
(sz^{s-1}+\Theta_{s}^{2})f_{nor}^{(m-1)}+(s\bar{z}^{s-1}+\Theta_{s}^{2})\overline{f_{nor}^{(m-1)}}
=\mathbb{P}_{N_s}^{m+1}.
$$
Hence
\begin{equation}
\begin{array}{lll}
g_{m+1}(w)&=&(1-s)^{(j+1)}\bar{a}z^{(j+1)s+1}(z\bar{z}+z^s)^{t-j-1}+(\bar{z}+sz^{s-1}
                +\Theta_{s}^{2})f_m(z,w)\\
          && +(z+s\bar{z}^{s-1}+\Theta_s^2)\overline{f_m(z,w)}+2Re\left((b_{N_0}-a_{N_0})z^{N_0}\right)+\Theta_{N_s}^{m+1}.
\end{array}
\end{equation}
By induction, we showed that if  the lemma holds for $m_0$ defined above, then it holds for any $m$ with $m_0 \leq m \leq 2t+(j+1)(s-2)+1$ and
$m\le N_0$.
\medskip

{\bf Step II of the proof of Lemma 2.3}:\  In this step, suppose that we know that the lemma holds  for   $m \in [
2t+j(s-2)+2,2t+(j+1)(s-2)+1]$ with $m\le N_0$, where $j$ is a certain non-negative integer bounded by $t-2$.
 We  then proceed to prove  that the lemma holds  also for $m
\in [ 2t+(j+1)(s-2)+2,2t+(j+2)(s-2)+1]$, whenever $m\le N_0$.

Suppose that $2t+(j+1)(s-2)+1<N_0$.
By the  assumption, we have
\begin{equation}
\begin{array}{lll}
g_{2t+(j+1)(s-2)+1}(w)&=&\bar{a}(1-s)^{j+1}z^{(j+1)s+1}(z\bar{z}+z^s)^{t-j-1}
     \\ && +(\bar{z}+sz^{s-1}+\Theta_{s}^{2})f_{2t+(j+1)(s-2)}(z,w)\\
     &&+(z+s\bar{z}^{s-1}+\Theta_s^2)\overline{f_{2t+(j+1)(s-2)}(z,w)}\\
   &&+2Re\left((b_{N_0}-a_{N_0})z^{N_0}\right)  +\Theta_{N_s}^{2t+(j+1)(s-2)+1}.
\label{eqn:Jam09}
\end{array}
\end{equation}
Collecting terms of degree $2t+(j+1)(s-2)+1$ in (\ref{eqn:Jam09}), we get
\begin{equation}
\begin{array}{lll}
g^{(2t+(j+1)(s-2)+1)}(z\bar{z})&=&\bar{a}(1-s)^{j+1}z^{(j+1)s+1}(z\bar{z})^{t-j-1}
+\hat{\mathbb{P}}_{N_s}^{2t+(j+1)(s-2)+1}\\
&&+\bar{z}f_{nor}^{(2t+(j+1)(s-2))}(z,z\bar{z})+z\overline{f_{nor}^{(2t+(j+1)(s-2))}(z,z\bar{z})}.
\label{eqn:Jam010}
\end{array}
\end{equation}
 Here we denote by $\hat{\mathbb{P}}_{N_s}^{2t+(j+1)(s-2)+1}$  a certain homogeneous
polynomial of degree $2t+(j+1)(s-2)+1$ with weight at least
$N_s$.

Now, we solve (\ref{eqn:Jam010}) as follows.
Write $\Lambda=2t+(j+1)(s-2)$.
Notice that
$$
I:=-\hat{\mathbb{P}}_{N_s}^{\Lambda+1}
+a(1-s)^{j+1}\bar{z}^{(j+1)s+1}(z\bar{z})^{t-j-1}+g^{(\Lambda+1)}(z\bar{z})
$$
is real valued and $I= \mathbb{P}_{N_s}^{\Lambda+1}$. Then (\ref{eqn:Jam010}) can be rewritten as
\begin{equation}
\begin{array}{lll}
I&=&\bar{a}(1-s)^{j+1}z^{(j+1)s+1}(z\bar{z})^{t-j-1}+a(1-s)^{j+1}\-{z}^{(j+1)s+1}(z\bar{z})^{t-j-1}
\\
&&+\bar{z}f_{nor}^{(2t+(j+1)(s-2))}(z,z\bar{z})+z\overline{f_{nor}^{(2t+(j+1)(s-2))}(z,z\bar{z})}.
\label{eqn:Jam010-2}
\end{array}
\end{equation}

Write
\begin{equation*}
I=\sum\limits_{\stackrel{ l+k=\Lambda+1}{l+(s-1)k \geq N_s}}
   a_{l\bar{k}}z^l\bar{z}^k.
\end{equation*}

Since $a_{l\bar{k}}=\overline{a_{k\bar{l}}}$, we also  require that $k+(s-1)l \geq N_s$.

We next can get the following general solution of (\ref{eqn:Jam010-2}):
\begin{equation}
\begin{array}{ll}
&f_{nor}^{(2t+(j+1)(s-2))}(z,w)=f^{(\Lambda)}_1+f^{(\Lambda)}_2\ \ \ \ \ \hbox{with}\\ & f^{(\Lambda)}_1=-\bar{a}(1-s)^{j+1}z^{(j+1)s+2}w^{t-j-2}
\\ & f^{(\Lambda)}_2=\sum_{\wt{l}+2\wt{k}=\Lambda}h_{\wt{l}\wt{k}}z^{\wt{l}}w^{\wt{k}}
\end{array}
\end{equation}
where $h_{\wt{l}\wt{k}}'s$ are determined by the following:
\begin{equation}
\sum\limits h_{\wt{l}\wt{k}}z^{\wt{l}+\wt{k}}\bar{z}^{\wt{k}+1}+\sum\limits
\overline{h_{\wt{l}\wt{k}}}z^{\wt{l}+\wt{k}+1}\bar{z}^{\wt{k}}=
\sum_{l,k}a_{{l}{k}}z^{{l}}\bar{z}^{{k}}.
\label{eqn:Jam12}
\end{equation}
Hence, we see that if $h_{\wt{l}\wt{k}}\not =0$, then either $\wt{l}=1, \ 2\wt{k}=\Lambda$ (in case $\Lambda$ is even)
or $\wt{l}+\wt{k}=l, \ \wt{k}+1=k.$ Here $l,\ k$ satisfy the properties described above.
Based on such an analysis and as argued before, we can conclude the following:
\begin{equation}
(sz^{s-1}+\Theta_{s}^{2})f_2^{(\Lambda)}(z,z\bar{z})
+(s\bar{z}^{s-1}+\Theta_s^2)\overline{f_2^{(\Lambda)}(z,z\bar{z})}
= {\Theta}_{N_s}^{\Lambda+2}.
\label{eqn:Jam13}
\end{equation}

Hence, from (\ref{eqn:Jam09})-(\ref{eqn:Jam13}), we get
\begin{equation}
\begin{array}{lll}
g_{\Lambda+2}(w)+g_{nor}^{(\Lambda+1)}(w)&=&(\bar{z}+sz^{s-1}+\Theta_{s}^{2})f_{\Lambda+1}(z,w)
            +(z+s\bar{z}^{s-1}+\Theta_s^2)\overline{f_{\Lambda+1}(z,w)}\\
       &&+\Theta_{N_s}^{\Lambda+2}+\hat{\mathbb{P}}_{N_s}^{\Lambda+1}
            +\bar{a}(1-s)^{j+1}z^{(j+1)s+1}(z\bar{z}+z^s)^{t-j-1}\\
       &&+(\bar{z}+sz^{s-1}+\Theta_{s}^{2})f_{nor}^{(\Lambda)}(z,w)
            +(z+s\bar{z}^{s-1}+\Theta_s^2)\overline{f_{nor}^{(\Lambda)}(z,w)}\\
       &&+2Re\left((b_{N_0}-a_{N_0})z^{N_0}\right).
\end{array}
\end{equation}
Notice that
\begin{equation*}
\begin{array}{ll}
&g_{nor}^{(\Lambda+1)}(z\bar{z})=\bar{a}(1-s)^{j+1}z^{(j+1)s+1}(z\bar{z})^{t-j-1}
    +\bar{z}f_{nor}^{(\Lambda)}(z,z\bar{z})+z\overline{f_{nor}^{(\Lambda)}(z,z\bar{z})}
    +\hat{\mathbb{P}}_{N_s}^{\Lambda+1},\\
&g_{nor}^{(\Lambda+1)}(w)-g_{nor}^{(\Lambda+1)}(z\bar{z}) \in
    \Theta_{N_s}^{\Lambda+2}.
\end{array}
\end{equation*}
We get\\
\begin{equation}
\begin{array}{ll}
g_{\Lambda+2}(w)=&(\bar{z}+sz^{s-1}+\Theta_{s}^{2})f_{\Lambda+1}(z,w)+2Re\left((b_{N_0}-a_{N_0})z^{N_0}\right)\\
     &+(z+s\bar{z}^{s-1}+\Theta_s^2)\overline{f_{\Lambda+1}(z,w)}
     +\Theta_{N_s}^{\Lambda+2}+J,
\end{array}
\end{equation}
where\\
\begin{equation}
\begin{array}{ll}
J=&(\bar{z}+sz^{s-1}+\Theta_{s}^{2})f_{nor}^{(\Lambda)}(z,w)
     +(z+s\bar{z}^{s-1}+\Theta_s^2)\overline{f_{nor}^{(\Lambda)}(z,w)}\\
  &
+\bar{a}(1-s)^{j+1}z^{(j+1)s+1}(z\bar{z}+z^s)^{t-j-1}-
     \bar{a}(1-s)^{j+1}z^{(j+1)s+1}(z\bar{z})^{t-j-1}\\
  &
-(\bar{z}f_{nor}^{(\Lambda)}(z,z\bar{z})+z\overline{f_{nor}^{(\Lambda)}(z,z\bar{z})}).
\end{array}
\end{equation}
Here we notice that
\begin{equation*}
\begin{array}{l}
\bar{z}f_{nor}^{(\Lambda)}(z,w)
     +z\overline{f_{nor}^{(\Lambda)}(z,w)}
-(\bar{z}f_{nor}^{(\Lambda)}(z,z\bar{z})+z\overline{f_{nor}^{(\Lambda)}(z,z\bar{z})})\\
\hskip 10pt+\bar{a}(1-s)^{j+1}z^{(j+1)s+1}(z\bar{z}+z^s)^{t-j-1}-
     \bar{a}(1-s)^{j+1}z^{(j+1)s+1}(z\bar{z})^{t-j-1}\\
=-\bar{a}(1-s)^{j+1}z^{(j+1)s+1}z\bar{z}(z\bar{z}+z^s)^{t-j-2}+\Theta_{N_s}^{\Lambda+2}\\
\hskip 10pt+\bar{a}(1-s)^{j+1}z^{(j+1)s+1}(z\bar{z})^{t-j-1}+\Theta_{N_s}^{\Lambda+2}\\
\hskip 10pt+\bar{a}(1-s)^{j+1}z^{(j+1)s+1}(z\bar{z}+z^s)^{t-j-1}-
     \bar{a}(1-s)^{j+1}z^{(j+1)s+1}(z\bar{z})^{t-j-1}\\
=\bar{a}(1-s)^{j+1}z^{(j+2)s+1}(z\bar{z}+z^s)^{t-j-2}+\Theta_{N_s}^{\Lambda+2}.
\end{array}
\end{equation*}
Hence we have

\begin{equation}
\begin{array}{lll}
J&=&
 (sz^{s-1}+\Theta_{s}^{2})f_1^{(\Lambda)}(z,w)
     +(s\bar{z}^{s-1}+\Theta_s^2)\overline{f_1^{\Lambda}(z,w)}
\\&&
+\bar{a}(1-s)^{j+1}z^{(j+2)s+1}(z\bar{z}+z^s)^{t-j-2}
+\Theta_{N_s}^{\Lambda+2}\\ &=&
\bar{a}(1-s)^{j+2}z^{(j+2)s+1}w^{t-j-2}+\Theta_{N_s}^{\Lambda+2}.
\label{eqn:Jam14}
\end{array}
\end{equation}
This proves the lemma when $m=2t+(j+2)s+2$. Now, the result obtained in the previous step completes the proof of the claim in this step.

\medskip
{\bf Step III of the proof of Lemma 2.3}:  We now can complete the
proof of the lemma by inductively using results obtained in Steps
I-II. Indeed, since we know that the Lemma holds for $m=2t+2$, we
see, by Step I, that the lemma holds for any $m\le N_0$ with $m\in
[2t+2, 2t+(s-2)+1]$.
 Then, applying first Step II and then  applying Step I again, we see the lemma holds for any $m\le N_0$ with $m\in [2t+j(s-2)+2, 2t+(j+1)(s-2)+1]$ and $j=1$. Now,  by
 an  induction argument on $j$, we see the proof of the lemma.
$\endpf$
\bigskip


We next complete the proof of Theorem 2.2 in case $\ord(f)=2t$.
First, if $m=ts+1<N_0$,  we then have, by Lemma 2.3:
\begin{eqnarray*}
\begin{array}{lll}
g_{ts+1}(w)&=&\bar{a}(1-s)^tz^{ts+1}+\Theta_{ts+2}^{ts+1}
  +(\bar{z}+sz^{s-1}+\Theta_{s}^{2})f_{ts}(z,w)\\
&&+(z+s\bar{z}^{s-1}+\Theta_s^2)\overline{f_{ts}(z,w)}.
\end{array}
\end{eqnarray*}
Collecting terms of degree $ts+1$ in the above equation, we obtain:

\begin{equation}
g_{nor}^{(ts+1)}(z\bar{z})=\bar{a}(1-s)^tz^{ts+1}
+\mathbb{P}_{ts+2}^{ts+1}+\bar{z}f_{nor}^{(ts)}(z,z\bar{z})+z\overline{f_{nor}^{(ts)}(z,z\bar{z})}.
\label{eqn:Jam15}
\end{equation}

Since $ts+2>ts+1$, we can write
$\mathbb{P}_{ts+2}^{ts+1}=\bar{z}A(z,\bar{z})$ for some polynomial function $A$. Hence, the
equation above is solvable only if a=0, which is a contradiction.

Second, suppose $2t+1 < N_0\le ts+1 $. By the normalization assumption in the theorem, we notice that $N_0\not = ts+1$. Hence,
we must have $2t+1 < N_0< ts+1 $

 Assume that $j$ is the integer such that
$ 2t+j(s-2)+2 \leq k_0s+j_0 \leq 2t+(j+1)(s-2)+1 $. Then by Lemma
2.3 and collecting terms of degree $N_0$ in $(\ref{eqn:Jam060})$, we
have
\begin{eqnarray*}
\begin{array}{lll}
g_{nor}^{(N_0)}(z\bar{z})&=&2Re\{(b_{N_0}-a_{N_0})z^{N_0}\}
  +\delta(1-s)^{j+1}\bar{a}z^{(j+1)s+1}(z\bar{z})^{t-j-1}
  \\ &&+\bar{z}f_{nor}^{(N_0)}(z,z\bar{z})+z\overline{f_{nor}^{(N_0)}(z,z\bar{z})}
  +\Theta_{N_0+1}^{N_0}.
\end{array}
\end{eqnarray*}
Here $\delta=0$ if $N_0 < 2t+(j+1)(s-2)+1$ and $\delta=1$
if
$N_0 = 2t+(j+1)(s-2)+1$.\\
With the same argument above, we can see a contradiction too.

Hence, to reach no contradiction, we must have $b_{N}=a_{N}$ for any $N\le ns+s-1$.
We thus conclude that $ts+1\ge ns+s$ and $t\ge  n+1$.
This  finally completes the proof.

\medskip
{\bf Step II of the proof of Theorem 2.2}:\ \ In this step, we show
that we can also have the result stated in Theorem 2.2 when
$\ord(f)$ is a finite odd number by applying the same argument as in
Step I.

Suppose that
 $\ord(f)=2t+1$,  then we can still assume that  $2t+2 \leq N_0$ as argued in Step I, where $N_0$ is defined  in a similar way.
Assume that $ts+s+1 < N_0$. Collecting  terms of degree 2t+2 in (\ref{eqn:Jam02}), we get
\begin{equation}
g_{nor}^{(2t+2)}(z\bar{z})=\bar{z}f_{nor}^{(2t+1)}(z,z\bar{z})+z\overline{f_{nor}^{(2t+1)}(z,z\bar{z})}.
\label{eqn:Jam20}
\end{equation}
Its solution is given by
\begin{equation}
f_{nor}^{(2t+1)}(z,w)=bzw^t \hskip 0.5cm , \hskip 0.5cm
g_{nor}^{(2t+2)}(w)=(b+\bar{b})w^{t+1}. \label{eqn:Jam21}
\end{equation}
Substituting the solution in (\ref{eqn:Jam21}) to  (\ref{eqn:Jam20})
and letting $A=(s-1)b-\bar{b}$, we get
\begin{equation}
\begin{array}{lll}
g_{2t+3}(w)&=&Az^s(z\bar{z}+z^s)^t+(\bar{z}+sz^{s-1}+\Theta_{s}^{2})f_{2t+2}(z,w)\\
            &&+(z+s\bar{z}^{s-1}+\Theta_s^2)\overline{f_{2t+2}(z,w)}
                +\Theta_{ts+s+1}^{2t+3}.
\label{eqn:Jam22}
\end{array}
\end{equation}
Repeating   the same induction argument as in the proof of Lemma 2.3, we get
\begin{equation}
\begin{array}{lll}
g_{ts+s}(w)&=&A(1-s)^tz^{ts+s}+\Theta_{ts+s+1}^{ts+s}+(\bar{z}+sz^{s-1}
                  +\Theta_{s}^{2})f_{ts+s-1}(z,w)\\
             &&+(z+s\bar{z}^{s-1}+\Theta_{s}^{2})\overline{f_{ts+s-1}(z,w)}.
\label{eqn:Jam22-1}
\end{array}
\end{equation}
Collecting terms of degree ts+s in (\ref{eqn:Jam22-1}), we obtain
\begin{equation}
g_{nor}^{(ts+s)}(z\bar{z})=A(1-s)^tz^{ts+s}+\mathbb{P}_{ts+s+1}^{ts+s}
+\bar{z}f_{nor}^{(ts+s-1)}(z,z\bar{z})+z\overline{f_{nor}^{(ts+s-1)}(z,z\bar{z})}.
\end{equation}
As before, it is solvable only when A=0 thus b=0, which gives a contradiction.
The case for $ts+s\ge N_0$ can be  similarly studied to conclude that $ts+s\ge ns+s$ and thus $2t+1\ge 2n+1$.
 This completes the proof of Theorem 2.2. $\endpf$

\bigskip

\begin{center}
{\section{\bf A complete set of formal invariants,  proofs of Theorem 1.1, Corollary 1.3  and  Theorem 1.5 }}
\end{center}

In this section, we will establish a formal normal form for the formal surface defined in (\ref{eqn:Jam000}), by applying a formal transformation
preserving the origin. This will give a complete classification of  germs of formal surfaces   $(M,0)$ with $\l=0,\ s<\infty$ in the formal setting, which, in particular, can  be  used to
answer an open question raised by J. Moser in 1985 ([pp 399, Mos]).

As another application of our complete set of formal invariants, we
show that a generic Bishop surface with the Bishop invariant
vanishing is not  equivalent to an algebraic surface, by
applying a Baire category argument similar to the study  in the CR
setting (see the paper of Forstneric [For]). Notice that this
phenomenon is strikingly different from the theory for elliptic
Bishop surfaces with non-vanishing Bishop invariants, where
Moser-Webster proved their celebrated theorem, that states that any
elliptic Bishop surface with a non-vanishing Bishop invariant has an
algebraic  normal form.


Let  $M$ be a formal Bishop surface in $\mathbb{C}^2$ defined by
\begin{equation}
w=H(z,\-{z})=z\bar{z}+2Re\{\sum_{j=s}^{N}a_jz^j\}+E_{N+1}(z,\bar{z}),
\label{eqn:Annie001}
\end{equation}
where $s\ge 3$ is a positive integer and   $E_{N+1}$ is a formal power series in $(z,\-{z})$  with $\ord(E_{N+1}) \geq N+1.$
Moreover, $a_s=1$ and for $ m>s,\ m\le N$,
$$a_m=0\ \ \  \hbox{if}\ \ m=0,1\ \hbox{mod s}.$$
 Our first  result of this section is the following normalization theorem:
\bigskip

{\bf Theorem 3.1}:  {\it With the above notation, there is a polynomial map
\begin{equation}
\left\{
\begin{array}{ll}
z'=z+f(z,w),\hskip 10pt&f(z,w)=O(|w|+|z|^2) \\
w'=w+g(z,w), &g(z,w)=O(|w|^2+|z|^3+|zw|)
\end{array}
\right.
\end{equation}
that transforms the  formal Bishop surface $M$ defined in (\ref{eqn:Annie001})
to the formal Bishop surface defined by
\begin{equation}
w'=H^*(z',\-{z'})=z'\bar{z'}+2Re\{\sum_{j=s}^{N+1}b_jz'^j\}+E^*_{N+2}(z',\bar{z'}).
\label {eqn:Annie003}
\end{equation}
Here $E^*_{N+2}=O(|z|^{N+2})$,  $a_j=b_j$ for $s\le j\le N$ and
$$b_{N+1}=0\ \  \hbox{if}\ \ N+1=0,1\ \hbox{mod s}.$$
Moreover, when $N+1\not =0,1\ \hbox{ mod}(s)$, $wt_{nor}(f)\ge N$ and $wt_{nor}(g)\ge N+1$;
when $N=ts$, $wt_{nor}(f)\ge 2t,\ \ wt_{nor}(g)\ge 2t+1.$
and when $N=ts-1$, $wt_{nor}(f)\ge 2t-1,\ \ wt_{nor}(g)\ge 2t.$}
\bigskip


Before proceeding to the proof, we recall a result of Moser, which will be used for our consideration here. For any $m\ge 4$ and holomorphic polynomials
$$f_{nor}^{(m-1)}(z,w),\ g_{nor}^{(m)}(z,w),\ \phi^{(m)}(z),$$ we define an operator, which we call  the Moser operator
${\mathcal L}$,  as follows:
$${\mathcal L}(f_{nor}^{(m-1)}(z,w), g_{nor}^{(m)}(z,w), \phi^{(m)}(z)):=g_{nor}^{(m)}(z,z\-{z})-2Re\{\-{z}f_{nor}^{(m-1)}(z,z\-{z})+\phi^{(m)}(z)\}.$$
The following lemma is essentially the content of Proposition 2.1 of
[Mos]:

\bigskip
{\bf Lemma 3.2}: Let $G(z,\-{z})$ be a homogeneous polynomial of degree $m$. Then
$${\mathcal L}(f_{nor}^{(m-1)}(z,w), g_{nor}^{(m)}(z,w), \phi^{(m)}(z))=G(z,\-{z})$$
has a unique solution: $\{f_{nor}^{(m-1)}(z,w), g_{nor}^{(m)}(z,w),
\phi^{(m)}\}$ under the normalization condition:
$f_{nor}^{(m-1)}=z^2f^*$ with $f^*$ a holomorphic polynomial.
Moreover, when $G$ has no harmonic terms, then ${\mathcal
L}\left(f_{nor}^{(m-1)}(z,w), g_{nor}^{(m)}(z,w),
0\right)=G(z,\-{z})$ also has a unique solution
$\{f_{nor}^{(m-1)}(z,w), g_{nor}^{(m)}(z,w)\}$ under the same
normalization condition just mentioned.
\bigskip

The proof of the Theorem 3.1 follows  from a similar induction argument that we used   in the previous section.

\bigskip
{\it Proof of Theorem 3.1}: We complete the proof in three steps.
\medskip

{\bf Step 1}: We first show that there is a polynomial map: $z'=z+f_{nor}^{(N)}(z,w),\ w'=w+g_{nor}^{(N+1)}(z,w)$, which maps
$M$ to a surface defined by the following equation:
\begin{equation}
w=z\bar{z}+2Re\{\sum_{j=s}^{N+1}b_jz^j\}+E^*_{N+2}(z,\bar{z})
\label{eqn:Annie004}
\end{equation}
 with $b_j=a_j$ for $s\le j\le N$.
Substituting the map into (\ref{eqn:Annie004}) and collecting terms of degree $N+1$, we see that the existence of the map is equivalent to the existence of solutions of the following functional equation:

\begin{equation}
{\mathcal L}(f_{nor}^{(N)}(z,w), g_{nor}^{(N+1)}(z,w), b_{N+1}z^{N+1})=-E^{(N+1)}_{N+1}(z,\-{z}).
\label{eqn:Annie005}
\end{equation}
By Lemma 3.2, we know that (\ref{eqn:Annie005}) is indeed solvable
and is uniquely solvable under the normalization condition as in
Lemma 3.2.

For the rest of the proof of the theorem, we can assume that $E_{N+1}=2Re\{b_{N+1}z^{N+1}\}+o(|z|^{N+1}).$
\medskip

{\bf Step 2}: In this step, we assume that $N+1=1\ \hbox{mod s}$. Write $N=ts$. We then show that there is a polynomial  map of the form:
\begin{equation}
\begin{array}{ll}
&z'=z+\sum_{l=0}^{N-2t}\{f^{(2t+l)}(z,w)\},\ \\ &w'=w+\sum_{\tau=0}^{N+1-2t-2}\{g_{nor}^{(2t+2+\tau)}(w)\}
\label {eqn:Annie006-02}
\end{array}
\end{equation}
such that under this transformation, $M$ is mapped to a formal surface $M'$ defined by (\ref{eqn:Annie003}) with $b_{N+1}=0$.
The map is also uniquely determined by imposing the normalization condition as in Lemma 3.2 for
$f^{(j)}$ with $2t<j\le N+1$.

As in Step I, this amounts to studying a series of normally weighted
homogeneous functional equations with the normally weighted degree
running from $2t$ to $N+1$. Substituting (\ref{eqn:Annie006-02})
into (\ref{eqn:Annie003}) and then collecting terms of degree
$2t+1$, we obtain the equation (\ref{eqn:Jam03}), which can be
solved as:
$$
f_{nor}^{(2t)}(z,w)=aw^t-\bar{a}z^2w^{t-1}$$
with $a$ to be (uniquely) determined later.

Now, suppose we are able to solve $f_{nor}^{(2t+l)},\
g_{nor}^{(2t+1+l)}$ for $2t+l=2t,\cdots, m-1\le st-2.$ Substituting
(\ref{eqn:Annie006-02}) into (\ref{eqn:Annie003}) and then
collecting terms of degree $m+1$, we  obtain an equation similar to
(\ref{eqn:Jam01-1}), as argued in the proof of Lemma 2.3:
\begin{equation}
g^{(m+1)}(z\bar{z})=\bar{z}f_{nor}^{(m)}(z,z\bar{z})+z\overline{f_{nor}^{(m)}(z,z\bar{z})}+\hat{\mathbb{P}}_{ts+2}^{m+1}
\label{eqn:Annie006-01}
\end{equation}
Notice that $\hat{\mathbb{P}}_{ts+2}^{m+1}(=\mathbb{P}_{ts+2}^{m+1})$ must be real valued and is uniquely determined by the known data.
This equation, in terms of the Moser operator, can be rewritten as:
\begin{equation}
{\mathcal L}\left(f_{nor}^{(m)}(z,z\bar{z}),g^{(m+1)}(z\bar{z}),0\right)=\hat{\mathbb{P}}_{ts+2}^{m+1}.
\label{eqn:Annie006}
\end{equation}
Since $\hat{\mathbb{P}}_{ts+2}^{m+1}$ is real-valued and divisible
by $\-{z}$, it   does not contain any harmonic terms. By Lemma 3.2,
it can be  solved,  and can be  uniquely solved under the
normalization condition in Lemma 3.2. By induction, we can uniquely
obtain $f_{nor}^{(m)},\ g_{nor}^{(m+1)}$  for $m\le ts-1$.
Substituting (\ref{eqn:Annie006-02}) into (\ref{eqn:Annie003}) and
then collecting terms of degree $m=ts+1$, we  obtain an equation
similar to (\ref{eqn:Jam15}), which can be rewritten as:
\begin{equation}
\begin{array}{ll}
&{\mathcal L}(g_{nor}^{(ts+1)}(z\bar{z}),f_{nor}^{(ts)}(z,z\-{z}),0) =2Re\{\bar{a}(1-s)^tz^{ts+1}\}
\\ &\hskip
1cm
+\hat{\mathbb{P}}_{ts+2}^{ts+1}-a(1-s)^t\-{z}^{ts+1}-2Re(b_{ts+1}z^{ts+1}).
\label{eqn:Annie007}
\end{array}
\end{equation}

As in the proof of Theorem 2.2,  the real-valued homogeneous
polynomial $\hat{\mathbb{P}_{ts+2}^{ts+1}}-a(1-s)^t\-{z}^{ts+1}$ has
a $\-{z}$ factor and thus has no harmonic terms. Hence, if we choose
$a=\-{b_{ts+1}}/(1-s)^t$, then (\ref{eqn:Annie007}) is uniquely solvable,
under the normalization condition in Lemma 3.2. This
completes the proof of the claim in this step.
\medskip

{\bf Step 3}:
 In this step, we assume that $N+1=0\ \hbox{mod s}$. Write $N=(t+1)s-1$.   We then show that  there is a unique  polynomial  map of the form:
\begin{equation}
\begin{array}{ll}
&z'=z+\sum_{l=0}^{N-1-2t}\{f_{nor}^{(2t+l+1)}(z,w)\},\ \\ &w'=w+\sum_{\tau=0}^{N+1-2t-2}\{g_{nor}^{(2t+2+\tau)}(w)\}
\label {eqn:Annie006}
\end{array}
\end{equation}
such that under this transformation, $M$ is mapped to a formal surface $M'$ defined by (\ref{eqn:Annie003}) with $b_{N+1}=0$. Here $f^{(m)}_{nor}$
satisfies the normalization condition in Lemma 3.2 for $m\not = 2t+1$.

The argument for this step is the same as  that for Step 2.
We first have to choose
$$f_{nor}^{(2t+1)}(z,w)=bzw^t, \hskip 1cm
g^{(2t+2)}_{nor}(w)=(b+\bar{b})w^{t+1}$$
with $b$ to be uniquely determined later. Arguing exactly in the same way as in Step 2, we can inductively find the unique solution (under the normalization condition) for
$f_{nor}^{(2t+l)},\ g_{nor}^{(2t+1+l)}$ with $2t+l=2t+2,\cdots, <st+s-1.$
At the level with  degree   $ts+s$, we have the following equation:

\begin{equation}
\begin{array}{ll}
&2Re(b_{N+1}z^{N+1})+g_{nor}^{(ts+s)}(z\bar{z})=((s-1)b-\bar{b})(1-s)^tz^{ts+s}
\\ &\hskip
1cm+\hat{\mathbb{P}}_{ts+s+1}^{ts+s}+\bar{z}f^{(ts+s-1)}(z,z\bar{z})
+z\overline{f^{(ts+s-1)}(z,z\bar{z})}. \label{eqn:Annie-002}
\end{array}
\end{equation}

Now, arguing the same way as in Step 2,  the equation
(\ref{eqn:Annie-002}) is uniquely solvable   by taking $b$ such that $(s-1)b-\-{b}=b_{N+1}$ and by imposing
the normalization condition as in Lemma 3.2  to $f^{(ts+s-1)}_{nor}$.

Now, the map in Theorem 3.1 can be chosen as  the map in Step 1 if $N+1\not = 0,1\ \h{mod}(s)$. When $N+1=0,\ \h{or}\ 1\ \h{mod}(s)$, the map in Theorem 3.1  can be defined  by composing  the map  in Step 2 or that in Step 3, respectively, with the map in Step 1.
We see  the proof of Theorem 3.1. Moreover, with such  fixed procedures and normalizations described in the above steps,
there are a set of universal polynomials $\{P_{kl}(a_{\a\b})\}_{1\le \a+\b\le k+l}$ (depending only on $s$ and $N$)
such that the coefficients of the map $(z',w')=(z,w)+(f,g)=(z,w)+\sum_{k,l}b_{kl}z^kw^l$  in Theorem 3.1 are determined by

\begin{equation}
b_{kl}=P_{kl}(a_{\a\b}),\ \ 1\le \a+\b\le k+l
\label{eqn:Annie-002-03}
\end{equation}
where $H=\sum_{\a,\b\ge 0}a_{\a\b}z^\a\-{z}^\b.$

The last sentence in Theorem 3.1 follows from the procedures that  we used to prove the existence part.
$\endpf$
\bigskip

We next choose the map $z'=z+f,\ w'=w+g$  in Theorem 3.1 such that its  coefficients are determined by (\ref{eqn:Annie-002-03}).
 Let $z=z'+f^*(z',w')$ and $w=w'+g^*(z',w')$ be its inverse transformation.
 Notice  that the coefficients of $(f^*,g^*)$ in its Taylor expansion up to degree, say $m$, are universal polynomial functions of the coefficients of
$(f,g)$ up to degree $m$ for any $m$. Hence we have the defining
equation of $M^*$, the image of $M$, as follows:
$$w'+g^*(z',w')=H(z'+f^*(z',w'), \-{z'+f^*(z',w')}).$$
Applying an implicit function theorem to solve for $w'$ and making
use of the uniqueness of the graph function, we see that the
coefficients in the Taylor expansion of $H^*$ up to degree $m$ must
also be polynomial functions of the coefficients of $H$ of degree
not exceeding $m$ in its Taylor expansion.
 Repeating such a normalization procedure that we did for $M$ to $M^*$ and by an induction argument, we get the following theorem: (The uniqueness part follows from Theorem 2.2.)
\bigskip

{\bf Theorem 3.3}: {\it Let $M$ be a formal Bishop surface defined by
\begin{equation}
w=H(z,\-{z})=z\bar{z}+z^s+\bar{z}^s+E(z,\bar{z}),
\label{eqn:Jam00}
\end{equation}
where $s\ge 3$ is a
positive integer and
$E(z,\-{z})=\sum_{\a+\b\ge s+1}^{\infty}a_{\a{\b}}z^\a\-{z^\b}$.
 Then there is a unique formal transformation of the form:
\begin{equation}
\left\{
\begin{array}{ll}
z'=z+f(z,w),\hskip 10pt&f(z,w)=O(|w|+|z|^2) \\
w'=w+g(z,w), &g(z,w)=O(|w|^2+|z|^3+|zw|)
\end{array}
\right.
\end{equation}
that transforms  $M$
to the formal Bishop surface defined by
\begin{equation}
w'=H^*(z',\-{z'})=z'\bar{z'}+z'^s+\-{z'}^s+ 2Re\{\sum_{ j=2,\cdots,s-1;\ k\ge 1}^{\infty}\lambda_{ks+j}z'^{ks+j}\}.
\label {eqn:Annie100}
\end{equation}
The normal form in (\ref{eqn:Annie100}),  up to a transformation of the form $z''=e^{i\theta}z',\ w''=w$  with
$e^{is\theta}=1$, uniquely determines the formal equivalence class of $M$.
Moreover,  there  are a set of universal polynomial functions
$$\{\Lambda_{ks+j}(Z_{\a\b})\}_{s+1\le \a+\b\le ks+j;\ j =2,\cdots,s-1;\ k\ge 1}$$
depending only on $s$, such that:
\begin{equation}
\begin{array}{ll}
\lambda_{ks+j}=\Lambda_{ks+j}(a_{\a\b})_{s+1\le \a+\b\le ks+j;\ j =2,\cdots,s-1;\ k\ge 1}.
\label{eqn:Annie101}
\end{array}
\end{equation}}
\bigskip

{\it Proofs of Theorem 1.1 and Corollary 1.3}:
Theorem 1.1 follows immediately from Theorem 3.3 and Lemma 2.1 (ii).

The proof of  Corollary 1.3 (a), (b), (d) also follows easily  from Theorem 3.1.
To see Corollary 1.3 (c), we let  $\mathcal G$ be a proper subgroup of ${\mathcal Z}_s$. Define $J_G:=\{j:\ 2\le j\le s-1,\ e^{i\theta j}=1,\
\h{for any}\ (e^{i\theta}z,w)\in {\mathcal G}\}$. Let $M_G$ be defined by
$$w=z\-{z}+z^s+\-{z}^s+2Re\{\sum_{j\in J_G}a_{s+j}z^{s+j}\},$$
with $a_{s+j}\not =0$.  Then we will  verify that $aut_0(M_G)={\mathcal G}$. To this aim, write
${\mathcal G}^*$ to be the collection of $\xi's$ with $(z,w)\ra (\xi z,w)$ belonging to ${\mathcal G}$.
By Corollary 1.3 (a), we need  only to show that if $\xi^{*s}=1$ and   $\xi^{*j}=1$ for any $j\in J_G$, then   $\xi^*\in {\mathcal G}^*$. Write $k=|{\mathcal G}^*|$. Then $s=km$ with $m(\in {\mathbf N})>1$. For any $\xi(\in {\mathcal G}^*)\not = 1$, since the order of $\xi$ must be divisible by $k$, we see that $\xi^k=1$. Therefore, ${\mathcal G}^*$ forms a complete set of the solutions
of $z^k=1$. Now, it is clear that $J_G=\{k,\cdots, (m-1)k\}$. Hence, we see that $\xi^{*k}=1$. Thus, $\xi^*\in {\mathcal G}^*$.
This completes the proof of Corollary 1.3 (c).

Now,  by Corollary 1.3 (a), we see that for $M$ as in Corollary 1.3
(e), $M$ must be formally equivalent to $M_s$. Assuming Theorem 1.2,
which we will prove in the next section, we also conclude that $M$
is biholomorphically equivalent to $M_s$. Corollary 1.3 (f) is a
simple consequence of the results in (a) and (e). $\endpf$
\bigskip

{\bf Corollary 3.4}: {\it Let $M$ be a real analytic Bishop surface defined by an equation of the form:
$$w=H(z,\-{z})=z\-{z}+2Re\{z^s+\sum_{k\ge 1,\ j=2,\cdots, s-1}a_{ks+j}z^{ks+j}\}\ \ \h{with infinitely many\ }\ a_{ks+j}\not = 0.$$
Then for any $N>s$, $M$ is not equivalent to the Bishop surface $M_N$ defined by
$$w=H_{(N+1)}(z,\-{z})=z\-{z}+2Re\{z^s+\sum^{ks+j\le N}_{k\ge 1,\ j=2,\cdots, s-1}a_{ks+j}z^{ks+j}\}.$$
 Here  $H_{(N+1)}$ is the $N^{th}$-truncation from the Taylor expansion of $H$ at $0$. In fact, $M_{(N+1)}$ is  equivalent to $M_{(N'+1)}$ with $N'>N$  if and only if   $a_{ks+j}= 0$  for any $N<ks+j\le N'$. }
\bigskip

Corollary 3.4 answers, in the negative,  the second  problem that J. Moser asked in his paper
([pp 399, Mos]).

\bigskip

As a less obvious application of Theorem 3.3, we next show that a
generic Bishop surface with the Bishop invariant vanishing at $0$
and with $s<\infty$ is not even formally equivalent to any algebraic
surface in ${\mathbb C}^2$. For this purpose, we borrow the idea
used in the CR setting based on the  Baire category argument. For
the consideration in the CR setting by using the Baire category
theorem, the reader is referred to the paper of Forstneric [For].

Write ${\mathcal M}_s$ for the collection of all formal Bishop surfaces defined as in (\ref{eqn:Jam00}):
\begin{equation}
w=H(z,\-{z})=z\bar{z}+2Re(z^s)+\sum_{\a+\b\ge s+1}a_{\a\b}z^\a\-{z^\b}.
\label{eqn:Annie001-1}
\end{equation}

 Write
${\mathcal F}:=\{\vec{a}=(a_1,\cdots,a_n,\cdots): \ a_j\in {\mathbb C}\}$,  equipped with the usual distance function:
$$dist(\vec{a},\vec{b})=\sum_{j=1}^{\infty}\frac{|a_j-b_j|}{2^j(1+|a_j-b_j|)}. $$
We know that $\mathcal F$ is a Fr\`echet space. There is a one-to-one correspondence between ${\mathcal M}_s$ and $\mathcal F$,
which assigns each $M\in {\mathcal M}_s$ to an element: $\vec{M}=(a_{\a\b})\in {\mathcal F}$  labeled in the lexicographical order. Therefore, we can, in what follows, identify
${\mathcal M}_s$  as a Fr\`echet space.
We define the operator $\mathcal J$ such  that it sends any $M\in {\mathcal M}_s$ to $(\lambda_{ks+j})_{j\not = 0,1;k\ge 1}$, where $(\l_{sk+j})$ is  described as in  Theorem 3.3. By (\ref{eqn:Annie101}), we easily see that $\mathcal J$ is a continuous map from ${\mathcal M}_s$ to $\mathcal F$.

 $(M,p)$ in ${\mathbb C}^2$ is called the germ of an algebraic surface  if $M$ near $p$ possesses a real polynomial defining equation. If $p\in M$ is a point with an elliptic  complex tangent, whose Bishop invariant is $0$ and whose Moser invariant is $s<\infty$, then there is a change of coordinates (see [Hu1], for instance) such that $p=0$ and
$M$ near $0$ is defined by an equation of the form:
\begin{equation}
w=z\-{z}+B(z,\-{z},w,\-{w}),\ \ B(z,\-{z},w,\-{w})=
\sum_{ 3\le \a+\b+2\gamma+2\tau}c_{\a\b\gamma\tau}z^\a \-{z}^\b w^{\gamma} \-{w}^\tau,
\label{eqn:Annie-00-0}
\end{equation}
where $B$ is a polynomial in its variables.
By using the implicit function theorem and using the argument in the step 1 of the proof of Theorem 3.1,
it is not hard to see that there is a fixed procedure to transform (\ref{eqn:Annie-00-0}) into a  surface defined by an equation  as in
 (\ref {eqn:Annie001-1}), in which $a_{\a\b}$ are  presented by polynomials of $c_{\a\b\gamma\tau}$ and $H(z,\-{z})$ becomes
what we call a Nash algebraic function to be defined as follows:

We call a real analytic function $h(z,\-{z})$ near $0$ a Nash
algebraic function if either $h\equiv 0$ or there is an irreducible
polynomial $P(z,\-{z};X)$ in $X$ with polynomial coefficients in
$(z,\-{z})$ such that $P(z,\-{z};h(z,\-{z}))\equiv 0.$ Certainly, we
can always assume that the coefficients of $(z,\xi,X)$ (in
$P(z,\xi,X)$) of terms with highest power in $X$ have  maximum value
$1$. The degree of $h$ is defined as the total degree of $P$ in
$(z,\-{z},X)$.

For $d,\ n,\ m\ge 1$, we define ${\mathcal A}^d_{B}(n,m)\subset {\mathcal M}_s$ to be the subset of Bishop surfaces defined in
(\ref{eqn:Annie001-1}), where $H(z,\-{z})'s$ are Nash algebraic functions derived  from the $B's$ in (\ref{eqn:Annie-00-0}) in the procedure described above with the degree of $B's$ bounded by $d$, that further satisfy  the following properties:

\medskip
\noindent
{\bf Cond (1)}: $H(z,\xi)'s$ are holomorphic over $|z|^2+|\xi|^2<1/m^2$;

\medskip

\noindent
{\bf Cond(2)}: $\max_{(|z|^2+|\xi|^2)<1/m^2}|H(z,\xi)|\le n$ and $|c_{\a\b\gamma\tau}|\le n$.
\medskip

Write ${\mathcal A}^d_{B}=\cup_{n,m=1}^{\infty}{\mathcal
A}^d_{B}(n,m)$ and ${\mathcal A}_B=\cup_{d=1}^{\infty}{\mathcal
A}^d_{B}$.
It is a consequence of Theorem 3.3 that
 $M$, defined in (\ref{eqn:Jam00}), is  formally equivalent to an algebraic surface if and only if ${\mathcal J}(M)\in {\mathcal J}({\mathcal A}_B)$.
(Therefore, $M$ defined in (\ref{eqn:Jam00}) is  not formally equivalent to an algebraic surface if and only if ${\mathcal J}(M)\not \in {\mathcal J}({\mathcal A}_B)$.)

Now, for any sequence $\{M_j\}\subset {\mathcal A}_B^d(n,m)$ with
$M_j: w=H_j(z,\-{z})=z\-{z}+z^s+\-{z}^s+o(|z|^s),$ by a normal
family argument and by passing to a subsequence, we can assume that
$H_j(z,\xi)\ra H_0(z,\-{z})$ over any compact subset of $\{|z|^2+
|\xi|^2<1/m^2\}$. If follows easily that $M_0$ defined by $w=H_0$ is
also in ${\mathcal A}_B^d(n,m)$. Moreover, $D^\a_zD^\b_\xi H_j(0)\ra
D^\a_zD^\b_\xi H_0(0)$ for any $(\a, \b)$. By (\ref{eqn:Annie101}),
${\mathcal J}(M_j)\ra {\mathcal J}(M_0)$ in the topology of
$\mathcal F$. Therefore, we easily see that ${\mathcal J}({\mathcal
A}_B)$ is a subset of $\mathcal F$ of the first category.

Next, for any $R>0$, we let
$${\mathcal S}_R:=\{\vec{\lambda}=(\l_{sk+j})_{k\ge 1;j = 2,\cdots,s-1}\}: \ \|\vec{\l}\|_R:=\sum_{ks+j}|\l_{ks+j}|R^{ks+j}<\infty\}.$$

It can be  verified that ${\mathcal S}_R$ is a Banach space under the above defined $\|\cdot\|_R$-norm. (In fact, it reduces to the standard $l^1$-space when $R=1$.)
We now claim that
${\mathcal K}^d_{B}$, defined as the closure of ${\mathcal J}\left({\mathcal A}^d_B(n,m)\right)\cap {\mathcal S}_R$ in ${\mathcal S}_R$ in its Banach norm, has no interior point.

Suppose, to the contrary,  that a  certain $\epsilon$-ball $\mathcal B$ of $\vec{a_0}=(\l^0_{sk+j})_{k\ge 1;j = 2,\cdots, s-1}$ in ${\mathcal S}_R$ is contained in
${\mathcal K}^d_{B}$. We must then have  ${\mathcal B}\subset {\mathcal J}\left({\mathcal A}^d_B(n,m)\right)\cap {\mathcal S}_R$. Indeed,
for any $\vec{a}\in {\mathcal B}$, let ${\mathcal J}(M_j)\ra \vec{a}$ with $M_j\in {\mathcal A}^d_{B}(n,m)$. By the argument in the above paragraph, we can assume, without loss of generality,  that
$M_j\ra M_0\in {\mathcal A}^d_{B}(n,m)$ in the $\mathcal F$-norm. By (\ref{eqn:Annie101}), we see that ${\mathcal J}(M_0)=\vec{a}$.
 Choose $\vec{a}=\{\l_{ks+j}\}$ such that $|\l_{ks+j}-\l^0_{ks+j}|\cdot(2R)^{ks+j}<\epsilon$ for any $ks+j$.  For any $N\ge 1$, then
we see that there is a certain $H=z\-{z}+z^s+\-{z^s}+\sum_{s+1\le \a+\b}a_{\a\b}z^\a\-{z^\b}$ Nash algebraic near $0$
such that
\begin{equation}
 \lambda_{ks+j}=\Lambda_{ks+j}(a_{\a\b}), \ \  N\ge ks+j\ge s+1, \ \a+\b\le ks+j,\ \ \Lambda=(\Lambda_{ks+j})_{s+1\le ks+j\le N}.
\label{eqn:Annie-002-02}
\end{equation}
Here $H$ is obtained from $B$ in (\ref{eqn:Annie-00-0}) with degree
of $B$ bounded by $d$. Since $a_{\a\b}$ are polynomial functions of
$c_{\a\b\gamma\tau}$, we can  conclude a contradiction from
(\ref{eqn:Annie-002-02}). Indeed,  since the variables on the right
hand side of (\ref{eqn:Annie-002-02}) are polynomially parametrized
by   less than $d^4$  free variables ($c_{\a\b\gamma\tau}$), the
image of  (\ref{eqn:Annie-002-02}) can not fill in an open subset of
${\mathbb R}^{N-s}$ as $N>>1.$


Therefore, we proved that ${\mathcal
A}_{B}=\cup_{d,n,m=1}^{\infty}{\mathcal A}^d_{B}(n,m)$ is a set of
the first category in ${\mathcal S}_R$. By the Baire category
theorem, we conclude  that most elements in  ${\mathcal S}_R$ are
not from  ${\mathcal J}\left({\mathcal A}_B\cap {\mathcal
S}_R\right)$. For any $\vec{a}=(\l_{sk+j})\not \in {\mathcal
J}\left({\mathcal A}_B\cap {\mathcal S}_R\right)$, the Bishop
surface defined by:
$w=z\-{z}+z^s+\-{z}^s+2\hbox{Re}(\sum_{k\ge1;j\not =
0,1}\l_{ks+j}z^{ks+j})$ is not equivalent to any algebraic surface
in ${\mathbb C}^2$.
 When $R$ varies,
we complete a proof of Theorem 1.5. $\endpf$
\bigskip

A real analytic  surface in ${\mathbb C}^2$ is called a Nash algebraic surface if it can be defined by a Nash algebraic function.
By the same token, we can similarly prove the following:

\bigskip
{\bf Theorem  3.5}: Most real analytic elliptic Bishop surfaces with
the Bishop invariant $\l=0$ and the Moser invariant $s<\infty$ at
$0$ are not equivalent to Nash algebraic surfaces in $\mathbb{C}^2$.
\bigskip

{\it Proof of Theorem 3.5}:
To prove Theorem 3.5, we  define ${\mathcal A}^d_{B}(n,m)$ in the same way as before except that we now only  require that
$H(z,\-{z})=z\-{z}+z^s+\-{z}^s+\sum_{\a+\b\ge s+1}a_{\a\b}z^{\a}\-{z}^\b$  is a general Nash algebraic function with total degree bounded by $d$ and with the same conditions described
as in Cond (1) and the first part of Cond (2). The last part of Cond (2) is replaced by the condition that $|b_{\a\b\gamma}|\le n$, where
$P(z,\-{z},X)=\sum b_{j}(z,\-{z})X^j=\sum_{\a\b\gamma}b_{\a\b\gamma}z^\a\-{z}^\b X^\gamma$ is a minimal polynomial of $H$ with the same coefficient restriction as imposed before.

We fix an $H_0$ and its minimal polynomial $P_0(z,\-{z};X)$. (We will fix certain coefficient of $P$ in the top degree terms of $X$ to be $1$ to make the minimal polynomial $P_0$ unique).
 Let  ${\mathcal A}^d_{B}(n,m; H_0,\delta)$ be a subset of ${\mathcal A}^d_{B}(n,m)$, where $M=\{ w=H(z,\-{z})\}\in {\mathcal A}^d_{B}(n,m; H_0,\delta)$ if and only if
 $|b_{\a\b\gamma}-b^0_{\a\b\gamma}|\le \delta$. Here $P=\sum b_{\a\b\gamma}z^\a\-{z}^\b X^\gamma$ and $P_0=\sum b^0_{\a\b\gamma}z^\a\-{z}^\b X^\gamma$ are the minimal polynomials of $H$ and $H_0$, respectively. We assume that $P$ is normalized in the same manner as for $P_0$. (Certainly, we can always do this if $\delta<<1$.)

Consider an $H$ and its minimal polynomial $P$ associated with an element from ${\mathcal A}^d_{B}(n,m; H_0,\delta)$.
Let $R$ be the resultant of $P$ and $P'_X$ with respect to $X$. We know that $R$ is a non-zero polynomial of $(z,\-{z})$ of degree bounded by $C_1(d)$, a constant depending only on $d$. Write $H=H^*_{(N)}+H^{**}_N$ with $H_{(N)}^{*}$ the Taylor polynomial
of $H$ up to  order $N-1$ and $H^{**}_N$ the remainder. Then from $P(z,\-{z},H^*_{(N)}+H^{**}_N)=0$, we obtain

\begin{equation}
P^{**}(z,\-{z}, X^{**})=0\ \ \hbox{ with}\ \ X^{**}=H_N^{**}.
\label{eqn:Annie-add}
\end{equation}

Here $P^{**}$ is a polynomial of total degree bounded by $C_2(d,
N)$, a constant depending only on $d$ and $N$, and its coefficients
are determined polynomially by the coefficients of $P$ and
$H_{(N)}^{*}$.  Notice that  $D_{X^{**}}\left(
P^{**}(z,\-{z},X^{**})\right)|_{X^{**}=0}=D_{X}\left (P(z,\-{z},
X)\right)_{X=H_{(N)}^{*}}$. Since there are polynomials $G_1$ and $G_2$
such that $G_1P+G_2P'_X=R$ and since
$P(z,\-{z},H^{*}_{(N)})=o(|z|^N)$, we conclude that the degree $k_0$
of the lowest non-vanishing order term of  $P'_X(z,\-{z},
H_{(N)}^{*})$ is  bounded by $C_1(d)$, depending only on $d$.

Choose an
$N\ge C_1(d)$ and a sufficiently small positive number $\delta$.
 We can apply a comparing coefficient method to (\ref{eqn:Annie-add}) to conclude that each $a_{\a_0\b_0}$ for $\a_0+\b_0\ge N$ is  determined by $b_{\a\b\gamma}$ and $a_{\a\b}$ with $\a+\b\le N-1$ through at most $C(k_0,N)$ rational functions in $b_{\a\b\gamma}$ and $a_{\a\b}$ ($\a+\b\le N-1$) with $C(k_0,N)$ depending only on $k_0, N$. Now,
(\ref{eqn:Annie-002-02}) can be used in the same manner to show that the interior of the closure of ${\mathcal J}({\mathcal A}^d_{B}(n,m;H_0,\delta))\cap {\mathcal S}_R$
in ${\mathcal S}_R$ is empty. It is easy to see that
 ${\mathcal J}({\mathcal A}_{B})$ can be written as  a countable union of these sets.
We see that ${\mathcal J}({\mathcal A}_{B})$
is a set of the first category in ${\mathcal S}_R$.
 This completes the proof of Theorem 3.5.
$\endpf$
\bigskip

{\bf Remark 3.6 (A)}:
The crucial point for  Theorem 3.5 to hold
is that the modular space of  surfaces with a vanishing Bishop invariant and $s<\infty$ is parameterized by an infinitely dimensional  space.
 Hence, any subclass of ${\mathcal M}_s$, that is  represented by
a countable union of finite dimensional subspaces of
${\mathcal M}_s$, is a thin set of ${\mathcal M}_s$
 under the equivalence relation. This  idea, that the infinite
 dimensionality of the modular space would generally have the
 consequence of  the generic non-algebraicity for its elements, dates
 back to  the early work of Poincar\'e [Po]. In  the CR setting,
 Forstneric in [Fo] has used the infinitely dimensional modular space
 of CR manifolds and the Baire category argument to give a short and quick
 proof that a generic CR submanifold in a complex space is not
 holomorphically equivalent to any algebraic manifold. Some earlier
 studies related to non-algebraicity for CR manifolds can be found,
 for instance, in
 [BER] [Hu2] [Ji].
However,  by a  result of the first author with Krantz [HK] and a result of the first
author in [Hu3], a Bishop surface with an elliptic complex tangent can
always be holomoirphically transformed into  the algebraic
Levi-flat hypersurface ${\mathbb C}\times {\mathbb R}$ and also into
the Heisenberg
hypersurface in ${\mathbb C}^2$.

{\bf (B)}. In the normal form (\ref{eqn:Annie100}), the condition that $\lambda_{ks+j}=0$ for $j=0,1,\ k=1,2,\cdots$ can be compared with the Cartan-Chern-Moser
chain condition in the case of strongly pseudoconvex hypersurfaces (see [CM]). In the hypersurface case, the chain condition
is  also  described by a finite system of differential equations. 
It would be very interesting to know if, in our setting here, there also exist a finite set of differential equations describing our chain condition. 

\bigskip
\begin{center}
{\section{\bf Surface hyperbolic geometry   and a convergence argument}}
\end{center}
In this section, we study the convergence problem for the formal
consideration in the previous section. Our starting point  is
the flattening theorem of Huang-Krantz [HK], which says that an
elliptic Bishop surface with a vanishing Bishop invariant can be
holomorphically mapped to ${\mathbb C}\times {\mathbb R}$.

Hence, to study the convergence problem, we can restrict ourselves to a real analytic  Bishop surface $M$ defined by
\begin{equation}
w=z\bar{z}+z^s+\bar{z}^s+E(z,\bar{z}),\ E(z,\bar{z})=\overline{E(z,\bar{z})}=o(|z|^s),\ 3\le s<\infty.
\label{eqn:Jenny-001}
\end{equation}
Recall that the Moser-Webster  complexification $\mathfrak{M}$ of $M$ is  the complex
surface near $0 \in \mathbb{C}^4$ defined by:\\
\begin{equation}
\left\{
\begin{array}{l}
w=z\zeta+z^s+\zeta^s+E(z,\zeta)\\
\eta=z\zeta+z^s+\zeta^s+E(z,\zeta).
\end{array}
\right.
\label{eqn:Jenny-002}
\end{equation}

We define the projection $\pi:\mathfrak{M}\longrightarrow
\mathbb{C}^2$ by sending $(z,\zeta,w,\eta) \in \mathfrak{M}$ to
$(z,w)$. Then $\pi$ is generically $s$ to 1. Write $B$ for the
branching locus of $\pi$. Namely, $(z,w) \in B$ if and only if
$\exists(\zeta_0,\eta_0)$ such that $(z,\zeta_0,w,\eta_0) \in
\mathfrak{M}$ and $\pi$ is not biholomorphic near
$(z,\zeta_0,w,\eta_0)$. Write ${\mathfrak B}=\pi^{-1}(B)$.

Then\\
\begin{equation*}
\begin{array}{l}
(z,w) \in B \\
\Longleftrightarrow \exists\zeta \ \hbox{ such  that }\
w=z\zeta+z^s+\zeta^s+E(z,\zeta)\ \hbox{ and }\
z+s\zeta^{s-1}+E_{\zeta}(z,\zeta)=0 \\
\Longleftrightarrow \sharp\{\pi^{-1}(z,w)\}<s.
\end{array}
\end{equation*}
It is easy to see that near $0,\ B$ is a holomorphic curve passing
through the origin.\\

Now, suppose $M'$ is defined by
$w'=z'\bar{z'}+z'^s+\bar{z'}^s+E^{\ast}(z',\bar{z'})$ with
$E^{\ast}(z',\bar{z'})=\overline{E^{\ast}(z',\bar{z'})}$ near 0.
Write $\mathfrak{M}'$ for the complexification of $M'$. Suppose that
$F:(M,0) \longrightarrow (M',0)$ is a biholomorphic map. Then $F$
induces a biholomorphic map $\mathcal {F}$ from $(\mathfrak{M},0)$
to $(\mathfrak{M}',0)$ such that $\pi' \circ \mathcal {F}=F \circ
\pi$. From this, it follows that $F(B)=B'$, where $B'$ is the
branching locus of $\pi'$ near the origin.

We next give the precise defining equation of $B$ near $0$.
From the equation $z+s\zeta^{s-1}+E_{\zeta}(z,\zeta)=0 $, we can
solve, by the implicit function theorem, that
\begin{equation}
z=h_1(\zeta)=-s\zeta^{s-1}+o(\zeta^{s-1}),
\label{eqn:Jenny-003}
\end{equation}
where $h_1(\zeta)$ is holomorphic near 0.
Substituting (\ref{eqn:Jenny-003}) into  (\ref{eqn:Jenny-002}), we get
\begin{equation}
w=h_2(\zeta)=(1-s)\zeta^s+o(\zeta^s).
\label{eqn:Jenny-004}
\end{equation}
From ( \ref{eqn:Jenny-004} ), we get
$$-\frac{w}{s-1}=(h_3(\zeta))^s\ \hbox{ with }\ h_3(\zeta)=\zeta+o(\zeta).$$
Hence, we get
\begin{equation}
\begin{array}{ll}
&\zeta=h_3^{-1}((-\frac{w}{s-1})^{\frac{1}{s}})
  =(-1)^{\frac{1}{s}}(\frac{1}{s-1})^{\frac{1}{s}}w^{\frac{1}{s}}+o(w^{\frac{1}{s}})\\
&z=h_1((-1)^{\frac{1}{s}}(\frac{w}{s-1})^{\frac{1}{s}}+o(w^{\frac{1}{s}}))=s(-1)^{-\frac{1}{s}}w^{\frac{s-1}{s}}\cdot(s-1)^{\frac{1-s}{s}}
  +o(w^{\frac{s-1}{s}}).
\end{array}
\end{equation}
Here, $h_j's$ are holomorphic functions near 0.
Next, let $w=u \geq 0$ and write
\begin{equation}
\begin{array}{ll}
&A_j(u)=h_1\circ h_3^{-1}\left (e^{-\frac{(2j+1)\pi
\sqrt{-1}}{s}}(\frac{u}{s-1})^{1/s}\right )\\
& =se^{\frac{(1+2j)\pi
\sqrt{-1}}{s}}u^{\frac{s-1}{s}}\cdot(s-1)^{\frac{1-s}{s}}
  +o(u^{\frac{s-1}{s}}),\ \ j=0,1,\cdots,s-1.
\end{array}
\label{eqn:Jenny-005}
\end{equation}

\textbf{Lemma 4.1}:\hspace{0.2cm} For $0 < u \ll 1,\ A_j(u) \in
D(u)$. Here
$$
D(u)=\{z \in \mathbb{C}^1: w=z\bar{z}+z^s+\bar{z}^s+E(z,\bar{z}) <
u\}.
$$
\bigskip

{\it Proof of Lemma 4.1}: The proof  follows clearly from the following estimate:
$$
|A_j(u)|^2+Re\{ 2A_j^s(u)+E(A_j(u),\overline{A_j(u)}) \}
=O(u^{\frac{2(s-1)}{s}}) \ll u
$$
as far as  $0<u \ll 1$ and $s \geq 3$.
$\endpf$
\bigskip

The following  fact will be crucial  for our later discussions:

\medskip
 $\{(A_j(u),u)\}_{j=0}^{s-1}=B \cap \{w=u\}$ and $A_j(u)$ is real analytic in $u^{1/s}$ for each fixed $j$.
\medskip

 Consider   a surface $(M, p)$ in ${\mathbb C}^2$. We say that $M$ near $p$ is defined   by a complex-valued function $\rho$,
 if $M$ near $p$ is precisely the zero set of $\rho$
 and  $\{Re(\rho),Im(\rho)\}$ has constant  rank two near $p$ as functions in $(x,y,u,v)$.
For a surface $(M,p)$ defined by $\rho$
 and a biholomorphic  map $F$ from a neigborhood of $p$ to a neighborhood of $p'$,
 we say that $F(M)$ approximates $(M^*,p')$ defined by $\rho^*=0$ to the order $m$  at $p'$ if there are smooth functions
$h_1$ and $h_2$ with $|h_1|^2-|h_2|^2\not = 0$ at $p'$
 such that $\rho\circ F^{-1}(Z)=h_1\cdot \rho^*+h_2\cdot\-{\rho^*}+o(|Z-p'|^m).$

\bigskip


\textbf{Lemma 4.2}:\hspace{0.2cm} Let $M,\ M'$ be  Bishop surfaces
near $0$ as defined above. Suppose that $F(M)$ approximates $M'$ to
the order $\wt{N}=Ns+s-1$ at $0$ with $N>1$. Then
$$
|F(A_j(u),u)-(A_j'(u'),u')| \ale |u|^{N}, \ \hbox{for} \
j=0,\cdots,s-1,\ u>0.
$$
Here $F=(z+f,w+g)$ is a holomorphic map with $f=O(|w|+|z|^2),\
g(z,w)=g(w)=O(w^2)$ and $u'=u+g(u)$.
\bigskip

{\it Proof of Lemma 4.2}:  Let $\Phi_1$ be a biholomorphic map, which maps M into
$M_{nor}^{N}$ defined by
$$
w=z\bar{z}+2Re\{z^s+\sum\limits_{k=1}^{N}\sum\limits_{j=2}^{s-1}a_{ks+j}z^{ks+j}\}
   +o(|z|^{sN+s-1}),
$$
and let $\Phi_2$ be a biholomorphic map from $M'$ to
${M'}_{nor}^{N}$
with ${M'}_{nor}^{N}$ defined by\\
$$
w'=z'\bar{z'}+2Re\{{z'}^s+\sum\limits_{k=1}^{N}\sum\limits_{j=2}^{s-1}{a'}_{ks+j}{z'}^{ks+j}\}
   +o(|z'|^{sN+s-1}).
$$
Define $\Psi=\Phi_2 \circ F \circ \Phi_1^{-1}$. Here we
assume
$\Phi_1,\ \Phi_2$  satisfy the normalization as in Theorem 3.1 at the origin.
Then $\Psi(M_{nor}^{N})$ approximates  ${M'}_{nor}^{N}$ up to order
$\wt{N}$.

By  Theorem 2.2, we conclude that
$$
a_{ks+j}=a'_{ks+j} \  for \  ks+j \leq {\wt{N}} \ \ \hbox{and}\  \
\Psi=Id+O(|(z,w)|^{ N}), \ \hbox{with}\ \ \wt{N}=Ns+s-1.
$$

In what follows, we write
$A_j(u),\ A_j^{\ast}(u),\ A_j^{nor}(u),\ A_j^{\ast nor}(u)$ for those
quantities, defined as in (\ref{eqn:Jenny-005}), corresponding to $M,\ M',\ M_{nor}^{N},\ {M'}_{nor}^{N}$,
respectively. Write $h_j^{nor}$ and $h_j^{\ast nor}$ for those holomorphic
functions, defined before, corresponding to $M_{nor}^{N}$ and ${M'}_{nor}^{N}$,
respectively. Then from the way these functions were constructed, we have
$$
h_j^{nor}(\zeta)=h_{j}^{\ast nor}(\zeta)+O(|\zeta|^{N}) \
for \ j=1,2,3.
$$
Hence,
$$
A_j^{nor}(u)=A_j^{\ast
nor}(u+\widetilde{g}(u))+O(u^{N}),
$$
where $\Psi=(z+\widetilde{f}(z),w+\widetilde{g}(w))$.
This immediately gives the following:
$$
F(A_j(u),u)=(A_j^{\ast}(u'),u')+O(u^{{N}}),
$$
where $u'=u+g(u),F=(z+f(z,w),w+g(w))$.
$\endpf$
\bigskip

Summarizing the above, we have the following:

\bigskip
\textbf{Proposition 4.3:}\hspace{0.2cm}{\it (1). Suppose that there
is a holomorphic map $F:M \longrightarrow M'$ with
$F=(z,w)+(O(|w|+|z|^2),O(w^2))$ such that $F(M)$ approximates $M'$
up to order $\wt{N}=Ns+s-1>s$ at 0. Then
$$
A_j^{\ast}(u+g(u))=A_j(u)+O(u^{N}), \ \ j=0,1,\cdots,s-1,\
u > 0.
$$

(2). Suppose that there is a formal holomorphic map $F:M
\longrightarrow M' $ with $F=(z,w)+(O(|w|+|z|^2),O(w^2))$. Then
$A_j^{\ast}(u+g(u))=f(A_j(u),u)$ in the formal sense. More
precisely, let $f_{(\wt{N})},\ g_{(\wt{N})}$ be the $(\wt{N}-1)^{th}$ truncation  in the
Taylor expansion of $f$ and $g$, respectively. Then
$$
A_j^{\ast}(u+g_{(\wt{N})}(u))-f_{(\wt{N})}(A_j(u),u)=O(u^{N'})
$$
where $N' \rightarrow \infty$ as $N \rightarrow \infty$. Indeed, we
can choose $N'=N$.}
\bigskip

Let $z=r\sigma(\tau,r)$ with  $u=r^2$ be the conformal map from the
 disk $r\Delta:=\{\tau\in{\mathbb C}: \ |\tau|<r\}$ to $D(u)$ with $\sigma(0,r)=0,\
\sigma_{\tau}'(0,r) > 0$. Here, as defined before, $$D(u)=\{z \in
\mathbb{C}^1: z\bar{z}+z^s+\bar{z}^s+E(z,\bar{z}) < u=r^2\}.$$
Similarly, let $z=r\sigma^*(\tau^*,r)$ with   $u=r^2$ be the
conformal map from the  disk $r\Delta$ to $D^*(u)$ with
$\sigma^*(0,r)=0,\ {\sigma^{*}}'_{\tau}(0,r) > 0$. Here,
$$D^*(u)=\{z \in \mathbb{C}^1: z\bar{z}+z^s+\bar{z}^s+E^*(z,\bar{z})
< u=r^2\}.$$ Then we know that $\sigma(\tau,r)=\tau(1+O(r))$ and
$\sigma$ is real analytic in $(\tau,r)$ over $\Delta_{1+\varepsilon}
\times (-\varepsilon,\varepsilon)$ with $0 < \varepsilon \ll 1$.
(See [Hu3]). Similar property also holds for $\sigma^*$.

Let $\tau_j(u) \in \Delta$ be such that
$r\sigma(\tau_j(u),r)=A_j(u)$. Then
$$
\tau_j(u)=\sigma^{-1}(\frac{A_j(u)}{u^{\frac{1}{2}}},\sqrt{u})
=\frac{A_j(u)}{u^{\frac{1}{2}}}(1+O(\sqrt{u}))
$$
Notice that
$\frac{A_j(u)}{u^{\frac{1}{2}}}=\sum\limits_{l=s-2}^{\infty}C_{l,j}u^{\frac{l}{2s}}$.
 Namely, $\frac{A_j(u)}{u^{\frac{1}{2}}}$ is analytic in
$u^{\frac{1}{2s}}$. Here
\begin{equation}
C_{s-2,j}=s(s-1)^{\frac{1-s}{s}}e^{\frac{\pi
\sqrt{-1}(1+2j)}{s}}. \label{eqn:Jenny-006}
\end{equation}

Now, the hyperbolic distance between $A_1(u)$ and $A_2(u)$ as points in
$D(u)$ is the same as the one between $\tau_1$ and $\tau_2$ as points in
$\Delta$.
Let $L_{1(j+1)}(u)=e^{d_{hyp}(\tau_0,\tau_j)}-1$.
In particular, $L_{12}(u)=e^{d_{hyp}(\tau_0,\tau_1)}-1$
Then
since
$$
  d_{hyp}(\tau_0,\tau_1)={1 \over
  2}\ln\left(\frac{1+|\frac{\tau_0-\tau_1}{1-\bar{\tau_0}\tau_1}|}
  {1-|\frac{\tau_0-\tau_1}{1-\bar{\tau_0}\tau_1}|}\right),\ \hbox{and}$$
$$
L_{12}(u)=s(s-1)^{\frac{1-s}{s}}
          |e^{\frac{\sqrt{-1}\pi}{s}}-e^{\frac{3\sqrt{-1}\pi}{s}}|
          u^{\frac{s-2}{2s}}+o(u^{\frac{s-2}{2s}}),
$$
we see that $L_{12}(u)$ is analytic in $u^{{1 \over 2s}}$.

Next, suppose $F:M \longrightarrow M'$ is a biholomorphic map with
$F=(\wt{f},\wt{g})=(z,w)+(O(|w|+|z|^2),O(w^2))$. Then
$\widetilde{f}=z+f$ is a conformal map from $D(u)$ to $D^*(u')$ with
$u'=u+g(u)$. Hence the hyperbolic distance between $A_1(u)$ to
$A_2(u)$ is the same as that of the hyperbolic distance from
$A_1^{\ast}(u')$ to $A_2^{\ast}(u')$, for
$F(A_j(u),u)=(A_j^{\ast}(u'),u')$.

Now, suppose that $F$ is a biholomorphic map with
$F=(\wt{f},\wt{g})=(z,w)+(O(|w|+|z|^2),O(w^2))$ such that $F(M)$
approximates $M'$ at 0 up to order $\wt{N}=Ns+s-1>s$. As before, we can assume
that $M,\ M'$ are already normalized up to order $\wt{N}$. Then
$F=Id+O(|z,w|^{N}),\
M=\{w=z\bar{z}+2Re\{\varphi_0(z)\}+o(|z|^{\wt{N}})\}$,
$M'=\{w=z\bar{z}+2Re\{\varphi_0(z)\}+o(|z|^{\wt{N}})\}$, where $\varphi_0
(z)=z^s+o(z^s),\ u'=u+g(u)=u+o(|u|^{N})$ and
$\varphi_0^{(sk+j)}(0)=0$ for $j = 0,1\ \hbox{mod  }(s).$

From the way $\sigma$ and $\sigma^*$ were constructed, we can show that (see [Lemma 2.1, Hu3]):
$$
\sigma^{\ast}(\tau,u')-\sigma(\tau,u)=\tau O(u^{N }).
$$
Indeed, this follows from the following more general result:
\bigskip

{\bf Lemma 4.4}: Let $\sigma(\xi,r)=\xi\cdot (1+O(r))$ and $\sigma^*(\xi,r)=\xi\cdot (1+O(r))$  be the biholomorphic map
from the unit disk $\Delta$ to
\begin{equation}
\begin{array}{ll}
&D(r):=\{\xi\in {\mathbb C}(\approx \-{\Delta}): \
|\xi|^2+rF_1(r,\xi,\-{\xi})<1\},\\ &D^*(r):=\{\xi\in {\mathbb C}
(\approx \-{\Delta}): \ |\xi|^2+rF_1(r,\xi,\-
{\xi})+r^mF_2(r,\xi,\-{\xi})<1\},
\end{array}
\end{equation}
respectively. Here $F_j(r,\xi,\-{\xi})$ are real-valued real
analytic functions in a neighborhood of $\{0\}\times \-{\D}
\times{\-\D}$. Then there is a constant $C$, depending only on
$F_j$, such that
$$
|\sigma^{\ast}(\xi,r)-\sigma(\xi,r)|\le C|\xi| r^{m},\ \ \xi\in\-{\D}.
$$

{\it Proof of Lemma 4.4}:  From the way $\sigma$ and $\sigma^*$ were constructed (see [Lemma 2.1, Hu3]), there
are $U, U^*\in C^{\omega}(\p\D\times (-\epsilon_0,\epsilon_0))$ with $0<\epsilon_0<<1$ such that
$$\sigma(\xi,r)=\xi\left(1+U(\xi,r)+{\mathcal H}(U(\cdot,r))\right),\ \sigma^*(\xi,r)=\xi\left(1+U^*(\xi,r)+{\mathcal H}(U^*(\cdot,r))\right), \ \xi\in \p\D.$$
Here ${\mathcal H}$ is the standard Hilbert transform and $U,\ U^*$ satisfy the following equations:
$$U=G_1(r,\xi,U,{\mathcal H} (U)),\ \ U^*=G_1(r,\xi,U^*,{\mathcal H}(U^*))+r^mG_2(r,\xi,U^*,{\mathcal H}(U^*)),$$
where $G_j(r,\xi,x,y)$ are real analytic in $(r,\xi,x,y)$ with $G_j \ale |r|+|x|^2+|y|^2$. Notice by the implicit function (see [Lemma
2.1, Hu3]), $\|U\|_{1/2},\ \|U^*\|_{1/2}\le C_1|r|$ with
$\|\cdot\|_{1/2}$ the H\"older-$\frac{1}{2}$ norm. Next, we have
\begin{equation}
\begin{array}{ll}
&U^*-U=\int_{0}^{1}\frac{\p G_1}{\p x}(r,\xi,\tau U^*+(1-\tau)U, \tau {\mathcal H} (U^*)+(1-\tau){\mathcal H}(U))(U^*-U)\hbox{d}\tau
\\ & +\int_{0}^{1}\frac{\p G_1}{\p y}(r,\xi,\tau U^*+(1-\tau)U,\tau {\mathcal H}(U^*)+(1-\tau){\mathcal H}( U))({\mathcal H}(U^*)-{\mathcal H}(U))\hbox{d}\tau
\\ & +r^mG_2(r,\xi,U^*,{\mathcal H}(U^*)).
\end{array}
\end{equation}
By noticing that the Hilbert transform is bounded acting on the H\"older space, we easily conclude the result in the lemma by letting $|r|<<1$.
$\endpf$
\bigskip

Now, by Lemma 4.4, we see that $\tau_j^{\ast}(u')=\tau_j(u)+O(u^{N})$.
Therefore
$$
L_{12}(u')-L_{12}(u)=O(u^{N }).
$$
In particular, if $F:M \longrightarrow M'$ is a formal equivalence
map
with $F=(\wt{f},\wt{g})=(z,w)+(O(|w|+|z|^2),O(w^2))$. Then
\begin{equation}
L_{12}^{\ast}(u')=L_{12}(u)\ \hbox{ in the formal sense}.
\label{eqn:Jenny-008}
\end{equation}

We next prove the following:

\bigskip

\textbf{Lemma 4.5:}\hspace{0.2cm} Let $F:M \longrightarrow M'$ be a
formal equivalence map with $F=(\wt{f},\wt{g})=(z,w)+(O(|w|+|z|^2),O(w^2))$. Write
$F=(\widetilde{f},\widetilde{g})=(z+f,w+g)$ as before. Then
$\widetilde{g}$ is convergent.
\bigskip

{\it Proof of Lemma 4.5}: By (\ref{eqn:Jenny-008}),
we have
$$
L_{12}^{\ast}(\widetilde{g}(u))=L_{12}(u)\ \ \hbox{in  the  formal
sense}.
$$
Write $u=V^{2s}$ and $\widetilde{g}(u)=U^{2s}$. Then
$$
L_{12}^{\ast}(U^{2s})=L_{12}(V^{2s}).
$$
Notice that $L_{12}^{\ast}(U^{2s})$ and $L_{12}(V^{2s})$ now are analytic
in $U$ and $V$, respectively. Moreover,
$$
L_{12}^{\ast}(U^{2s})=(\psi^{\ast}(U))^{s-2},\ L_{12}(V^{2s})=(\psi(V))^{s-2}
$$
with $\psi,\ \psi^{\ast}$ invertible holomorphic map of
$(\mathbb{C},0)$ to itself, and with $\psi'(0)={\psi^{\ast}}'(0)$.
Hence, we get
$$
\widetilde{g}(u)=((\psi^{\ast -1} \circ \psi)(u^{{1 \over
2s}}))^{2s}.
$$
On the other hand, $(\psi^{\ast -1} \circ
\psi(z)^{\frac{1}{2s}})^{2s}$ defines a multiple valued holomorphic
function near the origin. By the Puiseux expansion, we get
$$
(\psi^{\ast -1} \circ \psi(u^{{1 \over
2s}}))^{2s}=\sum\limits_{j=2s}^{\infty}c_ju^{{j \over 2s}}
$$
However, $(\psi^{\ast -1} \circ \psi(u^{{1 \over 2s}}))^{2s}$ also
admits a formal power series expansion. We conclude that $c_j=0$ if
$2s$ does not divide $j$. This proves the convergence of
$\widetilde{g}(u)$.
$\endpf$

\bigskip
We next prove the following theorem:
\bigskip

\textbf{Theorem 4.6:}\hspace{0.2cm} {\it Let $M,M'$ be  real analytic Bishop surfaces near $0$ defined by an equation of the form as in
(\ref{eqn:Jenny-001}). Suppose
$F=(\wt{f},\wt{g}):(M,0) \longrightarrow (M',0)$ is a formal equivalence map. Then F is biholomorphic
near $0$. }
\bigskip

{\it Proof of Theorem 4.6}:
We can  assume that $\wt{f}=z+f$ with $wt_{nor}(f)\ge 2$ and $\wt{g}=w+g$ with $wt_{nor}(g)\ge 4$.
By
 Lemma 4.4 and by considering $F_0\circ F$ instead of $F$, where $F_0(z,w)=(z,g^{-1}(w))$, we can assume, without loss of generality, that $\widetilde{g}=w$.
We will prove the convergence of $\widetilde{f}$ by the hyperbolic geometry associated to the surface discussed above.

By Proposition 4.3 (2), we first notice that
$$\widetilde{f}(A_j(u),u)=A_j^{\ast}(u)\ \hbox{ in the formal sense.}$$
 Namely,
$\widetilde{f_{(N)}}(A_j(u),u)=A_j^{\ast}(u)+o(u^{N'})$ for any $N$. Here,
$\widetilde{f_{(N)}}$ is the
 $N^{th}$-truncation in the Taylor expansion of $f$ at $0$;
and $N'$ depends only on N with $N' \rightarrow \infty$ as $N
\rightarrow \infty$.

Write $\widetilde{M}$ and $\widetilde{M'}$ for the holomorphic hull
of $M$ and $M'$, respectively. We next construct a holomorphic map
from $\widetilde{M} \setminus M$ to $\widetilde{M'} \setminus M'$ as
follows:

For $(z,u) \in D(u) \times \{u\}$, let $\tau(u) \in \Delta$ be such
that $r\sigma(\tau(u),r)=z,\ u=r^2$. Let $\Psi(\cdot,r)$ be a
biholomorphic  map from $\Delta$ to itself such that
$\Psi(\tau_j(u),r)=\tau_j^{\ast}(u)$ for $j=0,1$. Here, to see the
existence of $\Psi(\cdot,r)$ , it suffices for us to explain that
$d_{hyp}(\tau_0(u),\tau_1(u))=d_{hyp}(\tau_0^{\ast}(u),\tau_1^{\ast}(u))$.
But, this readily follows from (\ref{eqn:Jenny-008}) and Lemma 4.4.
Now, let
$$
\Psi_1=\frac{\tau-\tau_0(u)}{1-\bar{\tau_0}(u)\tau},
\Psi_1^{\ast}=\frac{\tau-\tau_0^{\ast}(u)}{1-\bar{\tau_0^{\ast}}(u)\tau},
\ \Theta(\tau,r)=e^{-i\theta(r)+i\theta^{\ast}(r)}\tau,
$$
where\\
\begin{equation*}
\begin{array}{l}
\theta(r)=
           arg\{\frac{\tau_1(u)-\tau_0(u)}{1-\overline{\tau_0}(u)\tau_1(u)}
           \frac{1}{u^{\frac{s-2}{2s}}}\} \\
\theta^{\ast}(r)=
           arg\{\frac{\tau_1^{\ast}(u)-\tau_0^{\ast}(u)}{1-\overline{\tau_0^{\ast}}
           (u)\tau_1^{\ast}(u)}\frac{1}{u^{\frac{s-2}{2s}}} \}.
\end{array}
\label{eqn:Jenny-010}
\end{equation*}
Then

\begin{equation}
\Psi(\tau,r)=\Psi_1^{\ast -1}(\tau,r) \circ \Theta(\tau,r) \circ
\Psi_1(\tau,r).
\label{eqn:Jenny-011}
\end{equation}

$\Psi(\tau,r)$ is analytic in $(\tau,u^{{1 \over 2s}})
\in \Delta_{1+\varepsilon_0} \times
(-\varepsilon_0,\varepsilon_0)$. (See [Lemma 2.1, Hu3]).

We notice that when $f$ is a priori known to be convergent, we then
have, by the uniqueness property of the M\"obius transformation,
that

\begin{equation}
\wt{f}(r\sigma(\xi,r),r^2)=r\sigma^{*}(\Psi(\xi,r),r^2).
\label{eqn:Jenny-012}
\end{equation}

Consider the angle $\Theta_j$ ($j=2,\cdots,s-1$)  from the geodesic connecting $\tau_j$ to $\tau_0$ to the geodesic connecting
 $\tau_j$ to $\tau_1$. As
a function of $u$ (or $r$), we see, as in the definition of $\Psi(\xi,r)$, that
$$\Theta_j(u)=arg\{\frac{\tau_1(u)-\tau_j(u)}{\tau_0(u)-\tau_j(u)}\cdot \frac{1-\-{\tau_j(u)}\tau_0(u)}{1-\-{\tau_j(u)}\tau_1(u)}\}=arg\{\frac{C_{s-2,2}-C_{s-2,j}}{C_{s-2,1}-C_{s-2,j}}\}+O(u^{1/(2s)}).$$
We can similarly define $\Theta_j^*$ for $M'$. Then the same argument which we used to show that $L_{12}(u)=L_{12}^*(u)$
can be used to prove that
$$\Theta_j(u)\equiv \Theta^*_j(u),\ \ \h{and} \ \ L_{1(j+1)}(u)=L^*_{1(j+1)}(u).$$

Now, we can use a M\"obius transformation to map $\tau_j$ to the
origin and $\tau_2$ to a point in the positive real line. Then we
easily see that $\Theta_j$ and $L_{1(j+1)}$ uniquely determine
$\tau_j(u)$. As an immediate consequence of such a consideration, we
conclude that

\bigskip

\textbf{Lemma 4.7}: \hspace{0.2cm}
$\Psi(\tau_j(u),r)=\tau_j^{\ast}(u) \ \ \hbox{for } j=0,\cdots,s-1$.
\bigskip

Now, for $(z,u) \in \widetilde{M}
\setminus M$ close to the
origin, we define
$$
f^{\ast}(z,u)=\sqrt{u}\sigma^{\ast}(\Psi(\sigma^{-1}(\frac{z}{\sqrt{u}},\sqrt{u}),\sqrt{u}),\sqrt{u})
$$
Then $f^{\ast}(z,u)$ is analytic in $\widetilde{M} \setminus M$.
We next prove the following:

\bigskip
\textbf{Lemma 4.8:}\hspace{0.2cm}$\forall \alpha \geq
0,\frac{\partial^{\alpha} f^{\ast}}{\partial
z^{\alpha}}(0,u)=\frac{\partial \widetilde{f}}{\partial
z^{\alpha}}(0,u)$ in the formal sense.  Namely, letting
$\widetilde{f_{(N)}}$ be the polynomial consisting of terms of
degree $\leq N$ in the Taylor expansion of $\widetilde{f}$ at $0$,
then $\exists N'(N)\rightarrow \infty$ as $N \rightarrow \infty$
such
that\\
$$
\frac{\partial^{\alpha} f^{\ast}}{\partial
z^{\alpha}}(0,u)=\frac{\partial^{\alpha} \widetilde{f_{(N)}}}{\partial
z^{\alpha}}(0,u)+o(u^{N'}).
$$
\bigskip

{\it Proof of  Lemma 4.8}: Let $S(u)$ be the hyperbolic polygon in $D(u)$
with vertices $A_j(u)(j=0,1,\cdots,s-1)$, whose sides consist of the geodesic segments connecting the vertices.
 Let $S^{\ast}(u)$ be
the one corresponding to $M'$. We notice that
 for any points $P,Q\in \D$,
then the geodesic segment connecting P to Q is
$$
\gamma_{P,Q}(t)=\frac{t\frac{Q-P}{1-Q\bar{P}}+P}{1+t\bar{P} \cdot
\frac{Q-P}{1-Q\bar{P}}}, \hskip 1cm 0\leq t \leq 1.
$$
Hence, by the same argument used in the proof of Lemma 4.2 and by making use of the property that $\wt{f}$ formally maps vertices to the corresponding ones,
 we see that for any point $P\in \p S(u)$,
we have
$$
f^{\ast}(P)
= \widetilde{f_{(N)}}(P)
+Error(P).$$

Here $$|Error(P)|\le C u^{N'} \ \ \hbox{with}\  \ N'(N)\rightarrow
\infty \ \hbox{as} \ N\rightarrow \infty
$$ and $C$ is a constant independent of $P$.

Now, by the Cauchy formula,
$$
\frac{\partial^{\alpha} f^{\ast}}{\partial z^{\alpha}}(0,u)=
\frac{\alpha!}{2\pi\sqrt{-1}}\int_{\partial
S(u)}\frac{f^{\ast}(\zeta,u)}{\zeta^{\alpha+1}}d\zeta
$$
and
$$
\frac{\partial^{\alpha} \widetilde{f_{(N)}}}{\partial z^{\alpha}}(0,u)=
\frac{\alpha!}{2\pi\sqrt{-1}}\int_{\partial
S(u)}\frac{\widetilde{f_{(N)}}(\zeta,u)}{\zeta^{\alpha+1}}d\zeta
$$
Notice that for $z \in \partial S(u),|z| \gtrsim u^{{s-1 \over s}}$
, it thus follows that
\begin{equation}
|\frac{\partial^{\alpha} f^{\ast}}{\partial z^{\alpha}}(0,u)-
\frac{\partial^{\alpha} \widetilde{f_{(N)}}}{\partial z^{\alpha}}(0,u) |
\ale O(u^{N'-{s-1 \over s}\alpha}).
\end{equation}
This completes the proof of Lemma 4.8. $\endpf$

\bigskip
We continue our proof of Theorem 4.6.
We notice that \\
$(i): \sigma^{\ast}(\zeta,\sqrt{u})$ is analytic in
$(\zeta,\sqrt{u})$ near $(0,0)$,\\
$(ii): \Psi(\tau,\sqrt{u})$ is analytic in
$\tau$ and $u^{{1 \over 2s}}$ near $(0,0)$ and,\\
$(iii): \sigma^{-1}(\frac{z}{\sqrt{u}},\sqrt{u})$ is analytic in
$(\frac{z}{\sqrt{u}},\sqrt{u})$ near $(0,0)$, too.\\
Write
$$
\Psi(\tau,\sqrt{u})=\sum\limits_{\alpha,\beta =
0}^{\infty}a_{\alpha\beta}\tau^{\alpha}u^{\frac{\beta}{2s}}
$$
and
$$
\Psi(\tau;Y_1)=\sum\limits_{\alpha,\beta =
0}^{\infty}a_{\alpha\beta}\tau^{\alpha}Y_1^{\beta}
$$
Then
$$
H(X,Y_1,Y_2)=Y_2\sigma^{\ast}(\Psi(\sigma^{-1}(X,Y_2);Y_1),Y_2)
$$
is analytic in $X,Y_1,Y_2$ near 0.
Write
\begin{equation}
H(X,Y_1,Y_2)=\sum\limits_{\alpha,\beta,\gamma=
0}^{\infty}b_{\alpha\beta\gamma}X^{\alpha}Y_1^{\beta}Y_2^{\gamma}
\label{eqn:Jenny-add-02}
\end{equation}
Then
$$
f^{\ast}(z,u)=H(\frac{z}{\sqrt{u}},u^{\frac{1}{2s}},\sqrt{u})
=\sum\limits_{\alpha,\beta,\gamma=
0}^{\infty}b_{\alpha\beta\gamma}z^{\alpha}u^{\frac{\gamma-\alpha}{2}
+\frac{\beta}{2s}}
$$
Hence,\\
$$
\frac{\partial^{\alpha} \widetilde{f}}{\partial z^{\alpha}}(0,u)=
\sum\limits_{\alpha,\beta,\gamma=
0}^{\infty}b_{\alpha\beta\gamma}\alpha! u^{\frac{\gamma-\alpha}{2}
+\frac{\beta}{2s}}
$$
in the formal sense.
It thus follows that if $b_{\alpha\beta\gamma} \neq
0,\frac{\gamma-\alpha}{2} +\frac{\beta}{2s}=\beta'$ is a
non-negative integer.
Hence $f^{\ast}(z,u)=\sum\limits_{\alpha,\beta,\gamma=
0}^{\infty}b_{\alpha\beta\gamma}z^{\alpha}u^{\beta'}$.\\
Now, $|b_{\alpha\beta\gamma}| \ale R^{|\alpha|+|\beta|+|\gamma|}$
with $R \gg 1$ by ( \ref{eqn:Jenny-add-02} ). We see that
$|b_{\alpha\beta\gamma}| \ale R^{2s|\alpha|+2s\beta'} \ale
(R^{2s})^{\alpha+\beta'}.$ This shows that $f^{\ast}(z,u)$ is
holomorphic in $(z,u)$ near 0; and thus $\widetilde{f}$ is
convergent, too. The proof of Theorem 4.6 is finally completed.
$\endpf$
\bigskip

{\it Proofs of Theorem 1.2 and Corollary 1.4}:  Theorem 1.2 and Theorem 4.6 have  the
same content.
The proof of Corollary 1.4  follows from Theorem 1.1 and Theorem 1.2.
$\endpf$

\bigskip

As another application of Theorem 1.2, we have the following
\bigskip

{\bf Corollary 4.9}: {\it Let $(M,0)$ be a real analytic elliptic
Bishop surface with the Bishop invariant vanishing  and the Moser
invariant $s<\infty$ at $0$. Then any element in $aut_0(M)$ is a
holomorphic automorphism of $(M,0)$. }


\bigskip

\vspace{10mm}

\begin{thebibliography}{widest-label}


\bibitem
[AG] {Gon 2006}
P. Ahern and X. Gong, Real analytic submanifolds in ${\mathbf C}^n$
with parabolic  complex tangents along a submanifold of codimension
one, preprint, 2006.



\bibitem[BER] {BER 2000}
S. Baouendi, P. Ebenfelt and L. Rothschild,
Local geometric properties of real submanifolds in complex space,
{\it Bull. Amer. Math. Soc.} (N.S.) 37 (2000), no. 3, 309--33.




\bibitem [BMR] {BMR 2002} S. Baouendi, N. Mir and  L. Rothschild,
Reflection ideals and mappings between generic submanifolds in
complex space, {\it J. Geom. Anal.  12}  (2002), 543-580.





\bibitem [BG] {BG 1984}
E. Bedford and B. Gaveau,
 Envelopes of holomorphy of certain
2-spheres in ${\mathbf C}^2$, {\it Amer. J. Math.} (105),  975-1009, 1983.


\bibitem [Bis] {Bis 1963} E. Bishop,
 Differentiable manifolds in complex
Euclidean space,  {\it Duke Math. J.} (32), 1-21, 1965.



\bibitem[Car]{Cartan 1932}
\'E. Cartan, Sur les vari\'et\'es pseudo-conformal des hypersurfaces de l'espace de deux variables complexes, {\it Ann. Mat.
Pura Appl. (4) 11}, 17-90(1932).




\bibitem [CM] {ChMo 1974}
S. S. Chern and J. K. Moser, Real hypersurfaces in complex manifolds,
{\it Acta Math. 133}, 219-271(1974).




\bibitem [Eli] {Eli 1997}
 Y. Eliashberg,
 Filling by holomorphic discs and its
applications, Geometry of Low-Dimensional Manifolds,
{\it London Math. Soc. Lecture
Notes Vol. 151},
 1997.









\bibitem [For] {For 2004} F. Forstneric,  Most real analytic
Cauchy-Riemann manifolds are nonalgebraizable, {\it Manuscripta Math. 115}
(2004),  489--494.





\bibitem
[Gon1] {Gon 1994} X. Gong,  On the convergence of normalizations of real analytic
 surfaces near hyperbolic complex tangents,
{\it Comment. Math. Helv}. 69 (1994), no. 4, 549--574.




\bibitem [Gon2] {Gon 1994-2}
X. Gong, Normal forms of real surfaces under
unimodular transformations near elliptic complex tangents, {\it Duke Math. J.
74}  (1994),  no. 1, 145--157.




\bibitem
[Gon3] {Gon 2004}
X.  Gong, Existence of real analytic
surfaces with hyperbolic complex tangent that are formally but not
holomorphically equivalent to quadrics, {\it Indiana Univ. Math. J. 53} (2004),
no. 1, 83--95.




\bibitem
[Grov] {Grov 1985}
 M. Gromov,
 Pseudo holomorphic curves in symplectic geometry, {\it Invent Math.
Vol. 82},
 1985,  307-347.






\bibitem [Hu1] {Hu 2004}
X. Huang, Local Equivalence Problems for Real Submanifolds in Complex Spaces,
{\it Lecture Notes in Mathematics 1848 (C.I.M.E. series)}, Springer-Verlag, pp 109-161, Berlin-Heidelberg-New York, 2004.


\bibitem
[Hu2] {Hu2001} X. Huang, On some problems in several complex
variables and Cauchy-Riemann Geometry, Proceedings of ICCM
(edited by L. Yang and S. T. Yau), {\it AMS/IP
Stud. Adv. Math. 20}, 383-396, 2001.


\bibitem
 [Hu3] {Hu3 1998} X. Huang, On  an n-manifold
 in ${\bf C}^n$ near an elliptic complex tangent,
{\it J. Amer. Math. Soc.} (11),  669--692,  1998.



\bibitem [HK] {HK 1995}
 X. Huang and S. Krantz,
 On a problem of Moser,
{\it Duke Math. J.} (78),
213-228, 1995.






\bibitem
[Ji] {ji} S. Ji,  Algebraicity of real analytic hypersurfaces with maxium rank, {\it Amer. Jour. Math. Vol 124}, 255-264, 2002.



\bibitem [KW1] {KW 1982}
C. Kenig and S. Webster,
 The local hull of holomorphy
of a surface in the space of two complex variables,
{\it Invent. Math.  67}, 1-21, 1982.


\bibitem [KW2] {KW2}
C. Kenig and S. Webster,
  On the hull of holomorphy of an n-manifold
in ${\mathbb C}^n$,
{\it Annali Scoula Norm. Sup. de Pisa IV  Vol. 11 (No. 2)},
 261-280, 1984


\bibitem [MMZ] {MMZ2003}
F. Meylan, N. Mir and D. Zaitsev, Approximation and convergence of formal CR-mappings, International Mathematics Research Notices 2003, no. 4, 211-242. 







\bibitem [Mos] {Mos 1985}
 J. Moser,
 Analytic surfaces in ${\CC}^2$ and their
local hull of holomorphy,
{\it Annales Aca -demi\ae Fennicae Series A.I. Mathematica}
(10), 397-410, 1985.


\bibitem [MW] {MW 1983}  J. Moser and S. Webster,
 Normal forms for real surfaces
in ${\CC}^2$ near complex tangents and hyperbolic
surface transformations,
{\it Acta Math.} (150), 255-296, 1983.


\bibitem [Po]
{Po} H. Poincar\'e, Les fonctions analytiques de deux variables et la
repr\'esentation conforme,
{\it Ren. Cire. Mat. Palermo, II. Ser. 23}, 185-220, 1907.


\bibitem [Sto]{Sto}
L. Stolovitch, Family of intersecting totally real manifolds of $(\Bbb C^n,0)$ and CR-singularities, preprint, 2006.


\bibitem [We] {Web}
 S. Webster, Pairs of Intersecting Real Manifolds in Complex Space, {\it Asian Jour. Math} Vol. 7 (No. 4), 449-462, 2003.

\bigskip
\noindent
X. Huang (huangx@math.rutgers.edu), School of Mathematics, Wuhan University, Wuhan 430072, China;
and Department of Mathematics, 
 Hill Center-Busch Campus, Rutgers University, 110 Frelinghuysen Road, Piscataway, NJ 08854-8019, USA;



\noindent
 Wanke Yin, School of Mathematical Sciences, Wuhan University, Wuhan 430072, China.

\end{thebibliography}
\end{document}